\documentclass[final]{siamltex}

\usepackage{graphicx}
\usepackage{float}
\usepackage{caption}
\usepackage{subfigure}
\usepackage{gensymb} 

\usepackage{amsfonts,amssymb}
\usepackage{color}
\definecolor{myblue}{rgb}{0.2, 0.5, 0.8}

\usepackage[pagewise]{lineno}
\usepackage{color}
\definecolor{darkgreen}{rgb}{0.1,0.5,0.1}
\definecolor{darkblue}{rgb}{0.2,0.2,1.0}
\usepackage[colorlinks=true,linkcolor=darkblue,citecolor=darkblue,
            filecolor=darkblue,urlcolor=darkgreen]{hyperref}

\usepackage{amsmath}
            
\newcommand{\p}{\partial}
\newcommand{\fr}{\frac}

\newenvironment{mat}{\left[ \begin{array}{ccccccccccccccc}}{\end{array}\right]}
\newenvironment{rmat}{\left[ \begin{array}{rrrrrrrrrrrrr}}{\end{array}\right]}
\def\bcm{\begin{mat}}
\def\ecm{\end{mat}}
\def\brm{\begin{rmat}}
\def\erm{\end{rmat}}
\newenvironment{pvect}{\left( \begin{array}{c}}{\end{array}\right)}
\def\bpvect{\begin{pvect}}
\def\epvect{\end{pvect}}

\newenvironment{choice}{\left\{ \begin{array}{ll}}{\end{array}\right.}

\def\choose{\begin{choice}}
\def\endch{\end{choice}}

\def\bsplit{\begin{split}}
\def\esplit{\end{split}}

\newcommand{\Fig}[1]{Figure~\ref{fig:#1}}

\newcommand{\ignore}[1]{}
\newcommand{\Comment}[1]{}

\newcommand{\eq}{\begin{equation}}
\newcommand{\en}{\end{equation}}
\newcommand{\eqm}{\begin{eqnarray}}
\newcommand{\enm}{\end{eqnarray}}
\newcommand{\eqmno}{\begin{eqnarray*}}
\newcommand{\enmno}{\end{eqnarray*}}
\newcommand{\eqml}[1]{\eql{#1}\begin{array}{rcl}}
\newcommand{\enml}{\end{array}\en}
\newcommand{\eql}{\begin{equation}\label}
\newcommand{\eqsub}[1]{\begin{subequations}\label{#1}\eqm }
\newcommand{\ensub}{\enm\end{subequations}}

\def\bc{\begin{center}}
\def\ec{\end{center}}
\def\bi{\begin{itemize}}
\def\ei{\end{itemize}}
\def\be{\begin{enumerate}}
\def\ee{\end{enumerate}}

\def\reals{{{\rm l} \kern -.15em {\rm R} }}

\def\qquad{\quad\quad}

\def\sf{{\cal F}}

\title{Computational and {\em in vitro} studies of
blast-induced blood-brain barrier disruption
\thanks{This work was supported in part by 
the Department of Veterans Affairs Office of
Research and Development Medical Research Service (DGC);
National Science and Technology Council of Mexico [CONACyT] (MdR),
NIH institutional fellowship [T32 AG000258] (JSM), Office of Academic
Affiliations, Advanced Fellowship Program in Mental Illness Research
and Treatment, Department of Veterans Affairs (BRH); NIH/NIA P50
AG005136 [ADRC] (ERP); VA RR\&D [I01RX001195] (ERP); University of
Washington Friends of Alzheimer's Research (DGC, ERP); 
NSF grant DMS-1216732 (RJL, MdR).}
}

\author{M.\ J.\ Del Razo\thanks{Department of Applied Mathematics, 
        University of Washington, Seattle, WA 98195-3925 ({\tt maojrs@uw.edu, rjl@uw.edu}).}
        \and Y.\ Morofuji \thanks{Department of Neurosurgery, University of Nagasaki, 
        Nagasaki, Japan ({\tt yoichi51@hotmail.com}).}
        \and J.\ S.\ Meabon \thanks{Northwest Network Mental Illness Research, Education, 
        and Clinical Center (MIRECC), VA Puget Sound Health Care System (VA Puget Sound), 
        Seattle, WA, USA and Department of Psychiatry and Behavioral Sciences, University of Washington 
        School of Medicine, Seattle, WA, USA ({\tt james64@u.washington.edu, bhuber@u.washington.edu, peskind@uw.edu}).}
        \and B.\ R.\ Huber \footnotemark[4]
        \and E.\ R.\ Peskind \footnotemark[4]
	\and W.\ A.\ Banks \thanks{Geriatric Research Education and Clinical Center (GRECC) 
	VA Puget Sound Health Care System (VA Puget Sound), Seattle, WA, USA 
	and Department of Medicine, University of Washington, Seattle, WA, USA ({\tt wabanks1@u.washington.edu, dgcook@uw.edu}).}
	\and P.\ D.\ Mourad \thanks{Department of Neurological Surgery, University of Washington, Seattle, WA USA
	and Division of Engineering and Mathematics, University of Washington, Seattle, WA, USA 
	({\tt doumitt@uw.edu}).}
        \and R.\ J.\ LeVeque \footnotemark[2]
        \and D.\ G.\ Cook \footnotemark[5]
        }

\begin{document}

\maketitle

\begin{abstract}
There is growing concern that blast-exposed individuals are at risk of
developing neurological disorders later in life. Therefore, it is important
to understand the dynamic properties of blast forces on brain cells,
including the endothelial cells that maintain the blood-brain
barrier (BBB), which regulates the
passage of nutrients into the brain and protects it from toxins in the blood.
To better understand the effect of shock waves on the BBB we
have investigated an {\em in vitro}
model in which BBB endothelial cells are grown in
transwell vessels and exposed in a shock tube, confirming
that BBB integrity is directly related to shock wave intensity.
It is difficult to directly measure the forces acting on these cells in the
transwell container during the experiments, and so a computational tool has
been developed and presented in this paper.

Two-dimensional axisymmetric Euler equations with the Tammann
equation of state were used to model the transwell materials, and a
high-resolution finite volume method based on Riemann solvers and the
Clawpack software was used to solve these equations in a mixed
Eulerian/Lagrangian frame. Results indicated that the geometry of the
transwell plays a significant role in the observed pressure time series in
these experiments.  We also found that
pressures can fall below vapor pressure due to the interaction of reflecting
and diffracting shock waves, suggesting that cavitation bubbles could be a
damage mechanism. Computations that include a simulated 
hydrophone inserted in the transwell suggest that the instrument itself
could significantly alter blast wave properties.
These findings illustrate
the need for further computational modeling studies aimed at understanding
possible blast-induced BBB damage.
\end{abstract}



\pagestyle{myheadings} \markboth{Computational and {\em in vitro} studies of
blast-induced BBB disruption}{M. J. Del Razo et al.}
\thispagestyle{plain}

\section{Introduction}\label{sec:intro}
Traumatic brain injury (TBI) is the leading cause of death and disability for 
people under the age of 45 years \cite{xiong2013animal}.
Non-penetrating impacts to the head are also associated with increased risk of
developing neurologic diseases that include Alzheimer's disease,
Parkinson's disease, and amyotrophic lateral sclerosis
\cite{fleminger2003head,plassman2000documented,bower2003head,lehman2012neurodegenerative}. 
In addition, repetitive mild traumatic brain injury (mTBI) has been
implicated in chronic traumatic encephalopathy 
\cite{mckee2014military,erlanger1999forum,mckee2009chronic,omalu2005chronic}. 
There is also growing evidence that repetitive low intensity non-impact blast wave exposure 
leads to mTBI, which similar to impact
TBI, can initiate slow-developing and potentially permanent brain disturbances \cite{cernak1997,cernak1996involvement,
cernak2001ultrastructural,mayorga1997pathology,murthy1979subdural,cernak2009traumatic,
warden2009case,huber2013blast,goldstein2012chronic}.  

The current and long-term health consequences of TBI and mTBI are of great
concern, particularly among military service members and Veterans, as well as 
civilian noncombatants \cite{al2008spectrum}.
Among US and coalition nations' military service members deployed to Iraq
and Afghanistan, it is estimated that approximately
15\% to 23\% have mTBI
\cite{wojcik2010traumatic,hoge2008mild,bell2009military,terrio2009traumatic}.
The majority of these mTBIs are blast-related \cite{bell2009military,
owens2008combat, peskind2011cerebrocerebellar, matthews2011multimodal,
dreier2014neuroimaging}, thus motivating the shock tube experiments described in
this paper.


Several sophisticated computational efforts (often employing commercial
finite element software) have been made in 
modeling TBI.
The majority of these efforts are aimed at
modeling the effects of a blast in an idealized human, mouse or rat head \cite{ruan1993finite,
anderson2003modern,takhounts2003development,takhounts2008investigation,kleiven2007predictors,
zhou1997viscoelastic,taylor2009simulation,wakeland2005computer,shreiber1997vivo,sundaramurthy2012blast},
sometimes including head-neck interactions \cite{kleiven2007predictors,ho2007dynamic,moss2009skull}. 
Much of this past work has been recently reviewed in \cite{gupta2013mathematical}. 

However, the mechanisms connecting blast wave exposure to mTBI are still not well understood. 
Clinical diagnostic neuroimaging approaches such as computerized tomography and
magnetic resonance imaging (MRI) fail to detect mild injuries.
This suggests that the injury mechanisms might occur at very small length
scales, even at the scale of a single cell. Several hypothesis have
been proposed: the disruption of BBB integrity \cite{shetty2014blood, hue2013blood, 
perez2013inflammatory}; cerebral vasospasm mechanotransduced
by the blast wave \cite{alford2011blast}; impairment of axonal functionality 
\cite{maia2014_2,maia2014_1}; shock wave excitation of phonons that decay into lower frequency 
oscillations \cite{kucherov2012blast} and the formation of cavitating bubbles
\cite{nyein2010silico,moss2009skull,moore2009computational,przekwas2009mathematical,
ziejewski2007selected,panzer2012development,goeller2012investigation}, among others.

In this paper, we study blast-induced BBB using differentiated brain-derived
microvessel endothelial cells, considered the biologically most relevant 
{\em in vitro}
approach for investigating BBB function \cite{hue2013blood,hue2014repeated}
(see Appendix for further discussion). 
Understanding such experiments is important since they isolate one possible
cause of TBI --- results 
show that blast functionally disturbs the BBB endothelial cell tight junction 
protein expression patterns. 
However, it is still extremely difficult to obtain accurate 
experimental measurements of the mechanical stresses exerted at the endothelial cells' location due 
to blast exposure, which could help relate specific damage mechanisms with experimental outcomes. 
In order to provide accurate quantitative data on the strength of the shock wave 
at this location, we developed a computational model which focuses on this 
particular experimental paradigm. 

The primary goal of this work is to computationally model the pressures to
which BBB endothelial cells
grown in a fluid-filled chamber placed in a shock tube are actually
subjected. The results obtained illustrate the fact that the geometry of the
chamber plays a large role in this, and suggest the possibility of
cavitation occurring in this experimental system.  More generally they
can aid in interpreting and understanding the experimental results. 
We also show that the introduction of a hydrophone into the experiment, as might be done
in an attempt to measure the pressures experimentally, could itself
change the outcome of the experiment and the likelihood of
cavitation occurring.

In an {\em in vivo} setting, the complexity of skull/bone 
anatomy, as well as the diffuse anatomy of the microvessel web in the brain, makes computational 
efforts to model BBB dysfunction extremely challenging. In contrast to this, the simple 
axisymmetric geometry of the 
{\em in vitro} system facilitates an accurate numerical investigation.
As explained further below in detail, this requires novel 
numerical algorithms to solve compressible Euler equations coupled with a Tammann equation of state (EOS)
across interfaces with large jumps in the material parameters at the interface between air and liquid. 
Numerical and exact methods for 
Euler equations with a Tammann EOS have been studied and developed previously, e.g.
\cite{Ivings1998,chertock2008interface,saurel1999simple,saurel2001multiphase,shyue1998efficient} 
among others; however, to the best of our knowledge, 
the numerical algorithms developed for this work are the only ones specifically developed to model
an experimental setup with fixed sharp interfaces with a big jump in the parameters. We present some description of the 
methods and a verification study. These methods are studied in more detail 
in \cite{delrazo2014,delrazo2015_02} and could be adapted to study related experiments.

\subsection{The biological effects of blast exposure on BBB
cells}
\label{sec:BBBexp}
One of the early manifestations of central nervous system (CNS) injury
following TBI is BBB disruption \cite{fujita2012intensity,
tomkins2011blood,rubovitch2011mouse}. The BBB is
responsible for maintaining and regulating separation between the CNS and
the circulating peripheral blood supply \cite{banks2010blood, zlokovic2011neurovascular}. 
In the brain, many cell types work together to regulate the BBB. However, the most 
important functional components of the BBB are the endothelial cells themselves, which
comprise the microvessels that supply the brain.  Brain endothelial cells
establish specialized connections called tight junctions with other
adjoining endothelial cells at points of cell-to-cell contact. This gives
rise to an extremely low-permeability cellular barrier that separates the
luminal (blood supply) side of the BBB from the abluminal (CNS) side of the
BBB. Significantly, there is evidence that BBB disruption may play an
important role in the delayed neurologic disorders associated with mTBI
\cite{shetty2014blood}. 

\begin{figure}
\hfill
\subfigure[Transwell]{\includegraphics[width=0.35\textwidth]{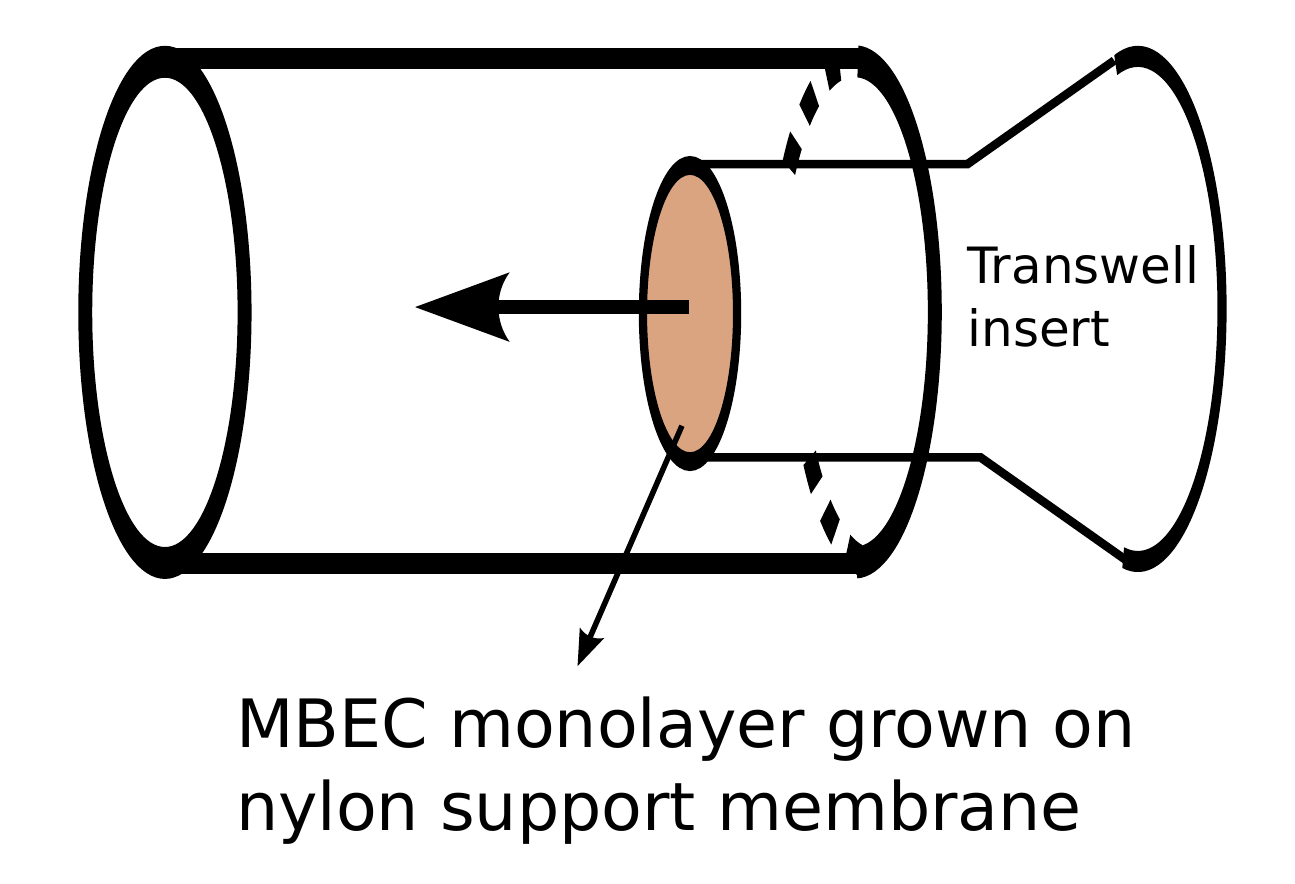}
\label{fig:stube-data-a}}
\hfill
\subfigure[Experimental setup]{\includegraphics[width=0.5\textwidth]{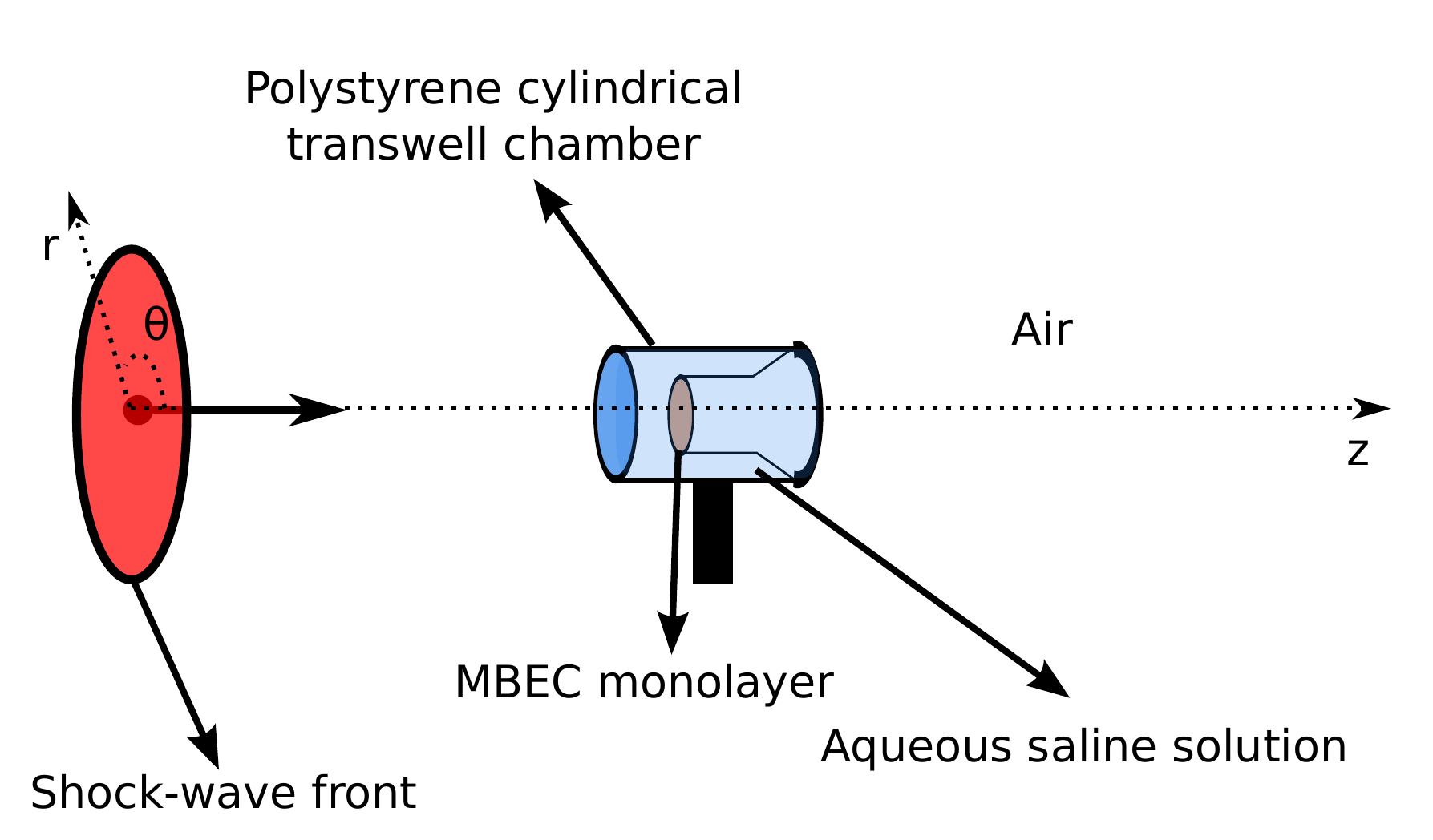}
\label{fig:stube-data-b}}
\vskip 5pt
\subfigure[Shock wave profile]{\includegraphics[width=0.6\textwidth]{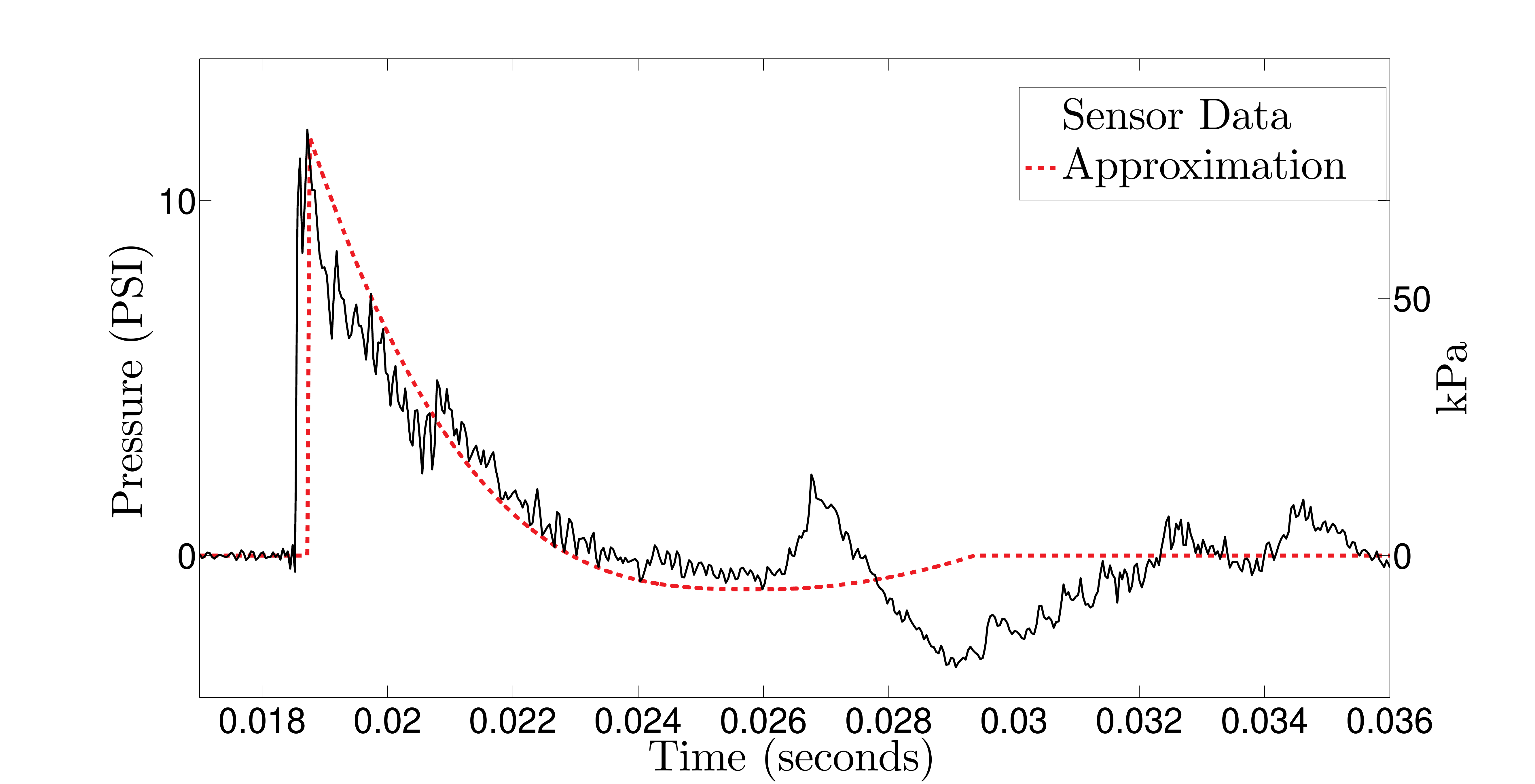}
\label{fig:stube-data-c}} 
\hfill
\subfigure[2D axisymmetric model]{\includegraphics[width=0.4\textwidth]{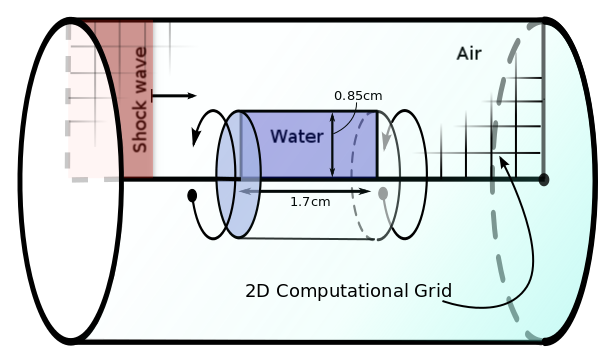}
\label{fig:stube-data-d}}
\hfill
\caption{(a) Polystyrene transwell chamber illustration. The transwell insert 
with the MBEC monolayers placed into the chamber filled with an aqueous solution.
(b) Cartoon of experimental system, showing the orientation of the transwell
in the shock tube.
The shock wave travels from the left through the air hitting the 
polystyrene transwell wall first, then the aqueous saline solution,
and finally the endothelial cells sample.
(c) The shock wave front profile obtained from a sensor 
before hitting the transwell as a function of time is shown as the
solid line. The approximation to be used as an initial
condition in the simulations herein is shown with a dashed line.
(d) The 3D axisymmetric shock tube model is obtained by revolving the
2D computational grid. The inside of the inner square corresponds to the cylindrical transwell
filled with aqueous saline solution, modeled here as water. The rest of the computational domain is a cylindrical
cross section of the shock tube filled with air.
}
\label{fig:stube-data}
\end{figure}

Recent studies have demonstrated that even mild blast exposures are capable
of disrupting the BBB \cite{abdul2013induction,yeoh2013distribution,li2013protective,
readnower2010increase}. In spite of this important progress, much
work remains in order to understand the mechanisms by which mild blast
exposure compromises BBB integrity. One approach to address this issue
is to study tight junctions using more simplified {\em in vitro models} of the BBB
\cite{banks2004triglycerides,banks2005effects,nakaoke2005human}. In this experiment, 
mouse brain-derived endothelial cells (MBECs) were isolated and
grown on permeable nylon support membranes, and then incubated in standard
cylindrical transwell tissue culture chambers, as illustrated in 
\Fig{stube-data-a}. Under these conditions, MBECs form an endothelial cell monolayer 
with mature tight junctions that functionally mimic the BBB \cite{banks2010blood,
zlokovic2011neurovascular}. The cylindrical transwell chamber was then 
completely filled with tissue culture media, sealed against leaks, placed 
inside a shock tube, and exposed to the blast, as shown in \Fig{stube-data-b}. 
Blast exposure has been shown to impair tight junction integrity under 
{\em in vitro} conditions, as well \cite{hue2013blood}.

Compared to {\em in vivo} conditions, in which the BBB is comprised of a highly 
elaborate matrix of microvessels in the brain, this {\em in vitro} BBB system offers a much
simpler geometry, with a planar MBEC monolayer positioned uniformly within a
defined cylindrical containment vessel (e.g. tissue culture chamber).

Although far removed from an actual brain, this {\em in vitro} approach provides 
the functional and anatomical precision required to correlate computed shock 
wave dynamics at a specific BBB that has a defined orientation with respect to 
propagating shock waves. 
Such combined anatomical and temporal precision is not possible under {\em in
vivo} experimental conditions. In addition, more complex computational models
of the brain cannot directly assess actual BBB biological function.
  
Importantly, the model presents new computational opportunities to better
estimate the biomechanical forces associated with blast overpressure exposure and thereby
derive more refined assessments of how forces elicited by blast exposure
affect BBB integrity under conditions that are biologically and independently
quantifiable. 

\begin{figure}[t]
\centering
\includegraphics[width=0.45\paperwidth]{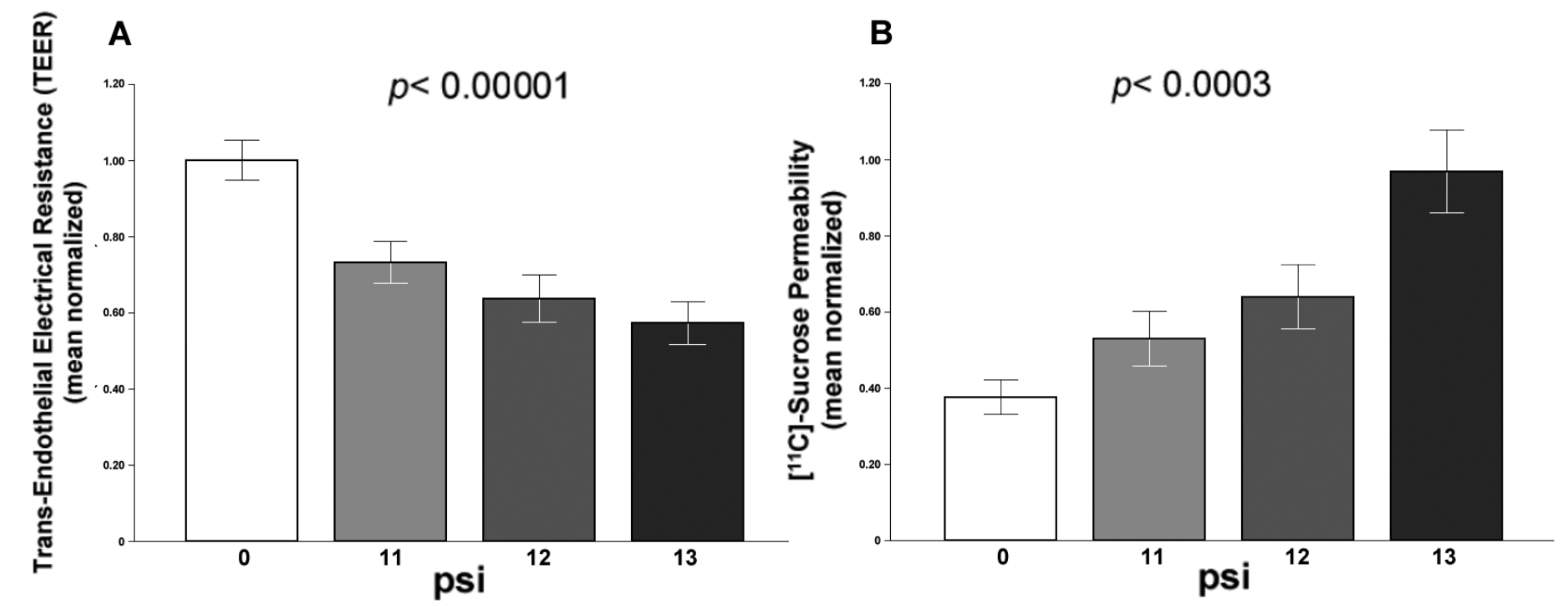}
\caption{
(A) Trans-endothelial electrical resistance (TEER) was significantly
decreased in a blast dose dependent fashion ($p \leq 0.00001$). Histogram
denotes mean normalized TEER at 24 hours after a sham exposure (0 psi)
($n$=26) or a single mild blast exposure with a peak amplitude of 11.0-11.9
psi (11), 12.0-12.9 psi (12), or 13.0-13.9 psi (13) ($n =$ 15, 12, and 9,
respectively). (B) MBEC monolayer permeability to radioactively labeled
[$^{14}$C]-sucrose was significantly increased in a blast dose dependent way
($p\leq 0.0003$) with the same blast exposure regimen as in panel A (blast 
intensity: 0, 11.0-11.9, 12.0-12.9, and 13.0-13.9 psi; $n$ = 7, 4, 3, and 4, 
respectively). Error bars indicate standard error of the mean (SEM).} 
\label{fig:BBBdiff}
\end{figure}

After exposure to the shock wave illustrated in \Fig{stube-data-b}, tests were
performed to measure the integrity of the BBB.
The results in \Fig{BBBdiff}A 
demonstrate that increasing blast intensity
produced a highly statistically significant decrease in trans-endothelial
electrical resistance (TEER) 24 hours post exposure 
($p\leq 0.00001$). In addition, there was a statistically significant negative
correlation between peak blast intensities (range: 0 -- 13.9 psi) and TEER(Pearson $r = -.603$, $p<0.00001$).

In a separate group of MBEC monolayers, we also measured blast-induced
leakage of [$^{14}$C]-labeled sucrose from the luminal transwell compartment
(i.e., peripheral circulating blood supply) into the abluminal transwell
compartment (i.e., CNS side). In keeping with the TEER measurements,
\Fig{BBBdiff}B shows that increasing blast intensity increased MBEC monolayer
permeability to [$^{11}$C]-sucrose ($p\leq 0.0003$). Consistent with this we found a statistically significant correlation between
overall peak blast intensities (range: 0--13.9 psi) and [$^{14}$C]-sucrose
permeability (Pearson $r = .695$, $p<0.001$). 

\subsection{The computational model}

The previous experimental results along with others presented in the Appendix
confirm that blast waves produce quantifiable and functional damage to BBB tissue. However, the
physical and/or biochemical mechanisms through which blast
damages brain tissue is not yet known. 
In order to gain insight on what some of these mechanisms might be, we have developed
a computational model based on the BBB experiment ---shown in \Fig{stube-data-b} and
described in the previous section--- that reproduces the dynamics and forces within the transwell 
chamber. The data computed with our model would be extremely difficult to obtain empirically, and moreover
the introduction of a measuring device would affect the outcome of the experiment,
as will be explored in detail in Section \ref{sec:compresults}. 

The computational model for this particular experiment consists of a
rectangular grid modeling the cylindrical 
axisymmetric cross-section of the shock tube.
A rectangular subsection of this grid models the polystyrene cylindrical transwell,
filled with saline solution (modeled as water), which is surrounded by air.  
The setup is shown in \Fig{stube-data-b}.

Some of the main issues that have been addressed with the computational
model presented in the next sections are:
\begin{itemize}
\item determine the shock wave interaction with an air-polystyrene-water interface, as in the experiment
from \Fig{stube-data-b}, to verify that the polystyrene layer can be omitted in
the computation;
\item explore the three-dimensional edge effects of the cylindrical transwell;
\item determine whether cavitation may be possible; 
\item explore how much the insertion of a hydrophone might modify measurements. \\
\end{itemize}

A necessary first step towards understanding the mechanical response of
BBB cells under shock loading is to 
determine the forces acting on the cells in the laboratory experiments. 
The shock strength 
increases as the shock passes from air into the fluid-filled transwell, but
the small diameter of the transwell results in waves also propagating in 
from the sides.  When the shock wave hits the distal end of the transwell a
reflected rarefaction wave is generated that interacts with the waves from
the sides and multiple wave reflections lead to a complex signal.

Moreover, the strong rarefaction waves propagating in the transwell
could result in fluid pressure values that are below the vapor
pressure, in which case cavitation bubbles may form.  
As cavitation bubbles collapse they can focus considerable kinetic energy
that is  capable of disrupting or destroying cellular membranes
\cite{ohl2006sonoporation,chen2011observations,xu2006new,
bigelow2009destruction,goeller2012investigation}.
Nonetheless, cavitation is a very complicated process that 
not only depends on the pressure but also on the amount of dissolved gas and other properties
of the fluid or tissue. Moreover, cavitation thresholds in the brain are variable 
and still largely unknown \cite{tung2010vivo,mcdannold2006targeted,tung2010identifying,mcdannold2008blood}. 
In this paper, we are only concerned with the possibility of cavitation in
the saline solution in the transwell, which was not de-gassed in the
experiment reported here.

The computational results obtained in the present paper --- although they 
provide limited answers --- support the possibility of
collapsing cavitation bubbles as one possible damage mechanism within the experimental 
arrangement. Note the algorithms and software developed are more widely applicable and could 
be adapted to study related experiments. For instance, 
cavitation could perhaps be directly modeled 
by extending these methods using a six-equation two-phase numerical 
model instead of Euler equations \cite{pelanti2014mixture,saurel2001multiphase}.

In the next section, we will show the results provided by the
computational model. In Section \ref{sec:nummethods}, we give
details of its numerical implementation, followed by further 
discussion in Section \ref{sec:disc}.

\section{Computational results} 
\label{sec:compresults} 
We will present the results of the computational version of the experiment shown in 
\Fig{stube-data-b}. The setup consists of the polystyrene transwell filled with saline solution,
modeled as water, without 
the endothelial cells, since these are too thin to be included in the model. 
Nonetheless, we can still measure the pressure intensity as a function of
time at the point where the cells are located. We will begin by citing a one-dimensional
version of the experiment done in a previous paper \cite{delrazo2014}, where we study the 
relevance of the thin polystyrene interface separating the air from the saline solution. Afterward,  
we will explore the full axisymmetric two-dimensional model that 
will allow us to study the edge effects and possible cavitation. Finally, we repeat this 
experiment with the addition of an hydrophone-shaped inclusion in order to
determine how the inclusion of such a pressure-measuring device might affect the experiment. 

\subsection{Air-polystyrene-water interface}

In a previous work \cite{delrazo2014}, we implemented a one-dimensional 
version of the experimental system in \Fig{stube-data-b} by zooming in on
the left face of the transwell chamber. The one-dimensional model consists of 
only three interfaces: air, polystyrene and water. Since the
polystyrene walls of the transwell are very thin relative to the characteristic 
length of the experiment, we study the effect of decreasing the thickness of the 
polystyrene layer on the transmitted shock wave. We show that when the polystyrene 
interface is thin enough in comparison to the transwell length, the 
results are effectively the same as without it. This result allows us to set up our 
two-dimensional axisymmetric model with only one fixed interface between air and 
saline solution and completely neglect the effect of the polystyrene walls. The results and 
methods from this section are explained in more detail elsewhere \cite{delrazo2014}.

\subsection{Two-dimensional axisymmetric results: Cavitation and edge effects}
With these simplifications in mind, we constructed the two-dimensional axisymmetric computational model. 
The implementation was done using the methods of Section \ref{sec:nummethods} 
to solve the two-dimensional axisymmetric Euler equations (\ref{eq:Eulercyl}) coupled with
the Tammann equation of state (EOS) (\ref{eq:SGEOS}) to model the different materials. The three-dimensional 
solution is recovered from revolving the solution on the two-dimensional grid as shown in \Fig{stube-data-d}, so the 
model is effectively three-dimensional. The geometry of the air and water interfaces is also shown. 
The air and water parameters for the Tammann EOS 
are the ones given in Table \ref{tab:param}. Furthermore, to provide an accurate model, we need to model 
length scales according to the experiment. The cylindrical transwell 
filled with water (saline solution) is $1.7$ cm long with a radius of $0.85$cm; it can be modeled as a two-dimensional rectangle before
being revolved. The shock wave is modeled by feeding the profile shown in \Fig{stube-data-c} to the
left boundary of the computational domain. However, on the time and length scales of the simulation, we only observe the 
shock wave and an essentially constant pressure behind the shock, since the rarefaction wave that reduces the pressure
behind the shock wave decays over roughly 3 msec while the computation is run for only 134 $\mu$s.

\begin{figure}[t]
  \centering
  \includegraphics[width=0.32\textwidth]{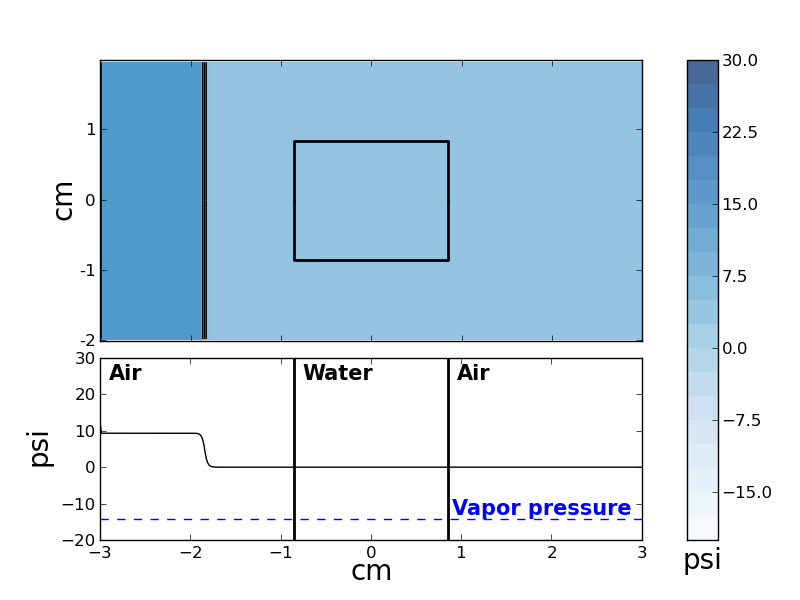} 
        \includegraphics[width=0.32\textwidth]{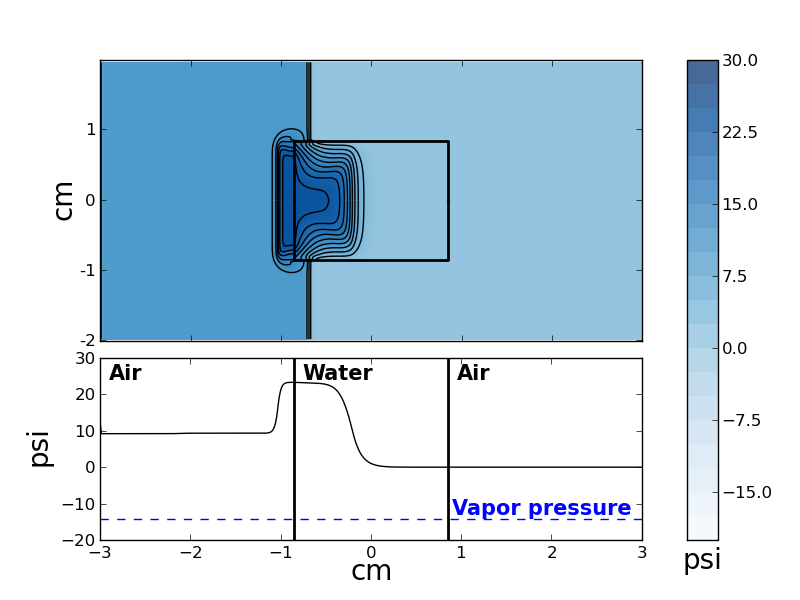} 
        \includegraphics[width=0.32\textwidth]{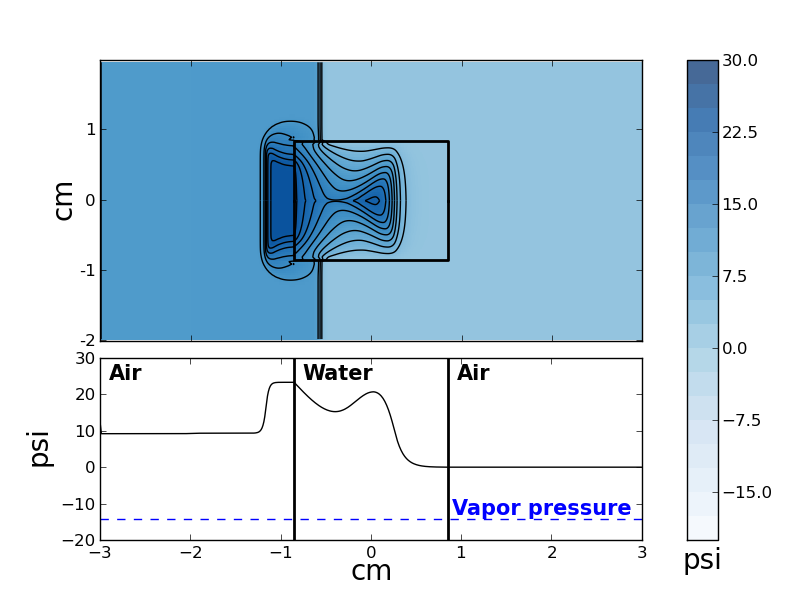} 

  \includegraphics[width=0.32\textwidth]{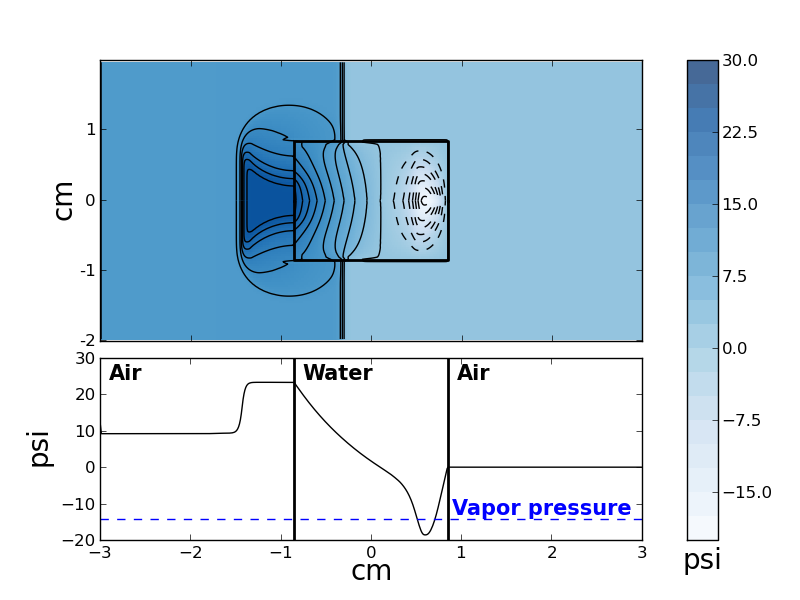} 
        \includegraphics[width=0.32\textwidth]{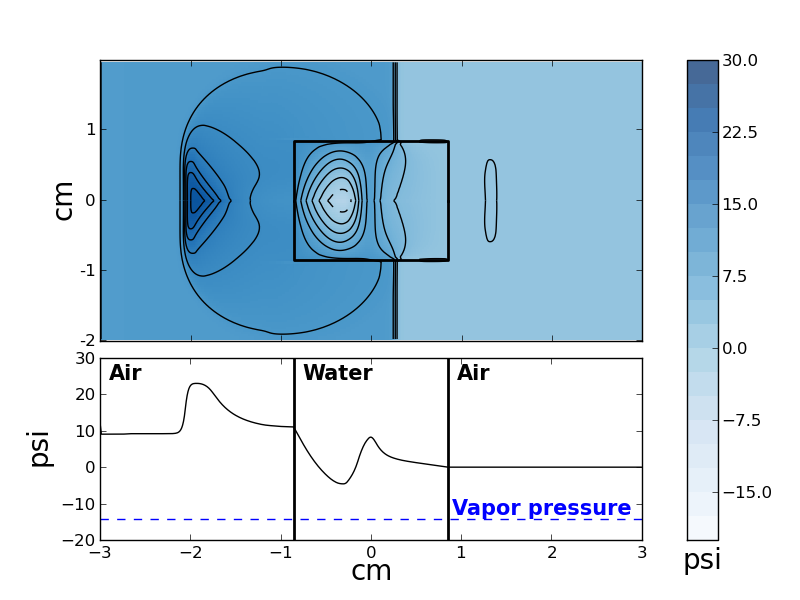} 
        \includegraphics[width=0.32\textwidth]{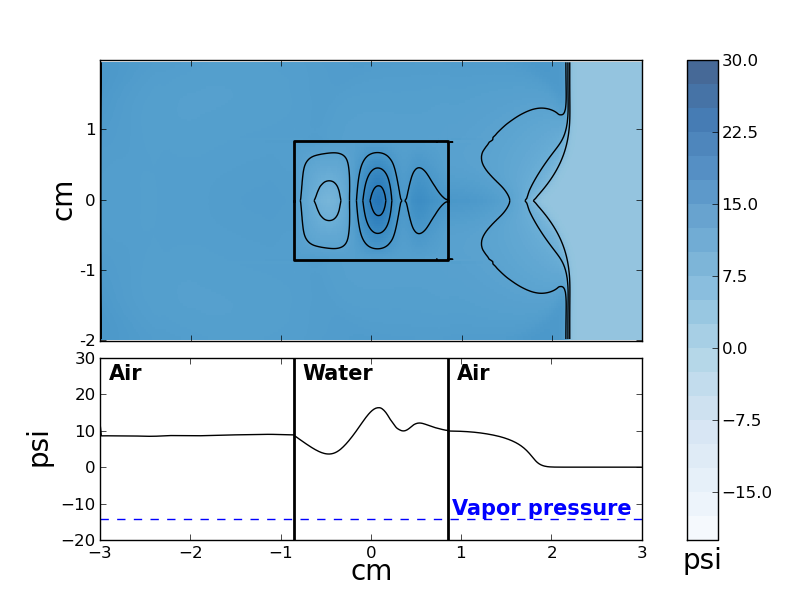} 
  \caption{Axisymmetric simulation output at six different times points
  $t=30,60,63.2,69.6,84.8,134.4$ $\mu s$. Two-dimensional pressure contour plots of a planar 
  cross section of the cylinder are shown, along with pressure trace along the axis.
  Water vapor pressure is also shown indicating where cavitation 
  might be possible. Distance is displayed in centimeters and pressure in psi, where atmospheric pressure 
  corresponds to 0 psi. }
  \label{fig:TBIanim3D}

  \includegraphics[width=0.32\textwidth]{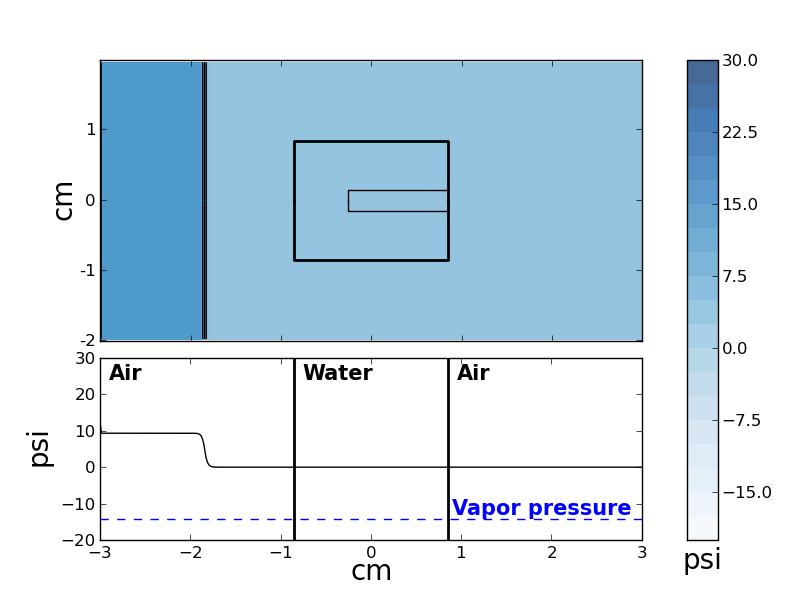} 
        \includegraphics[width=0.32\textwidth]{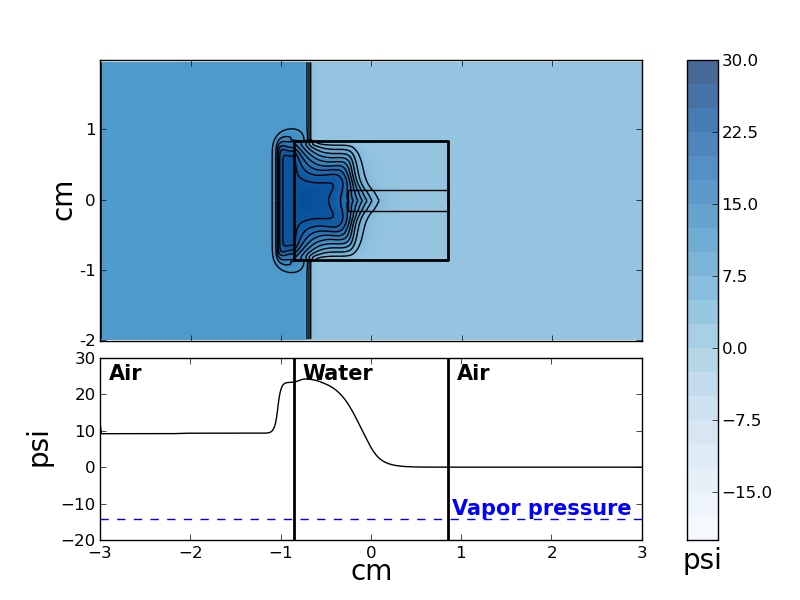} 
        \includegraphics[width=0.32\textwidth]{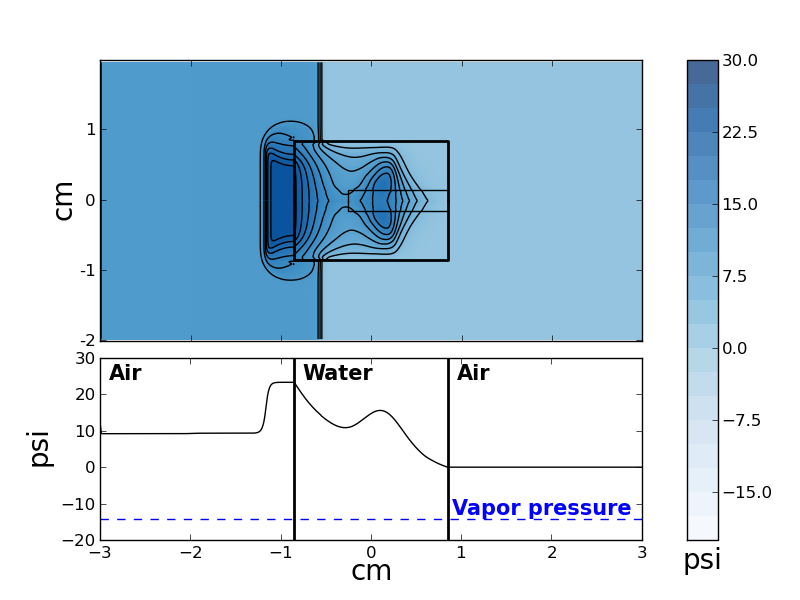} 
          \label{fig:TBIanim2D1_hyd}
  \includegraphics[width=0.32\textwidth]{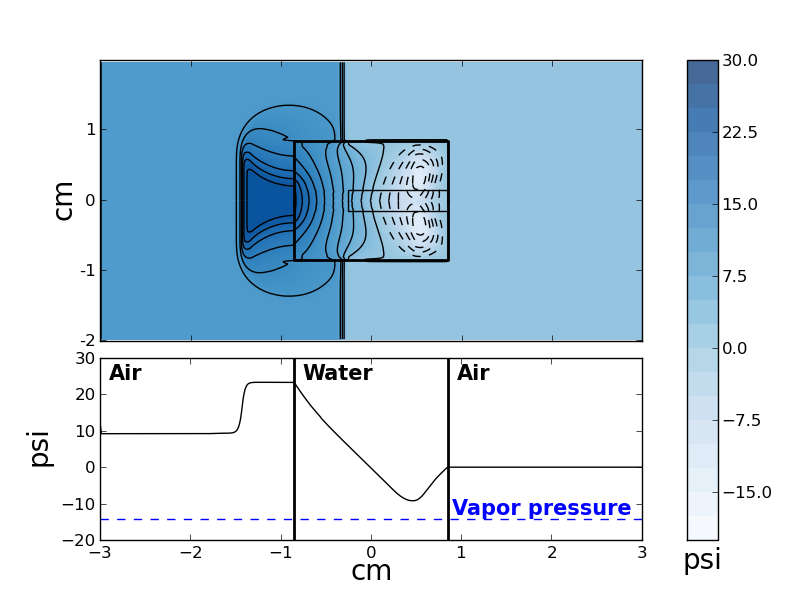} 
        \includegraphics[width=0.32\textwidth]{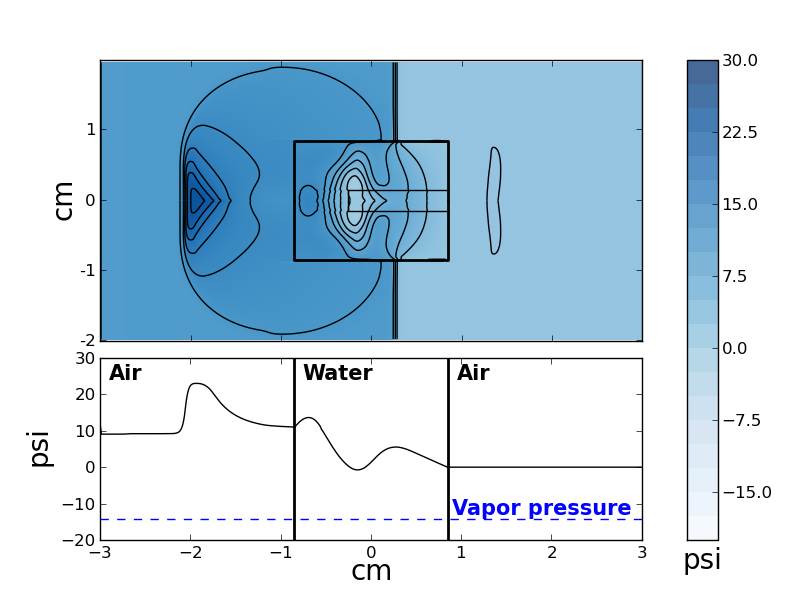} 
        \includegraphics[width=0.32\textwidth]{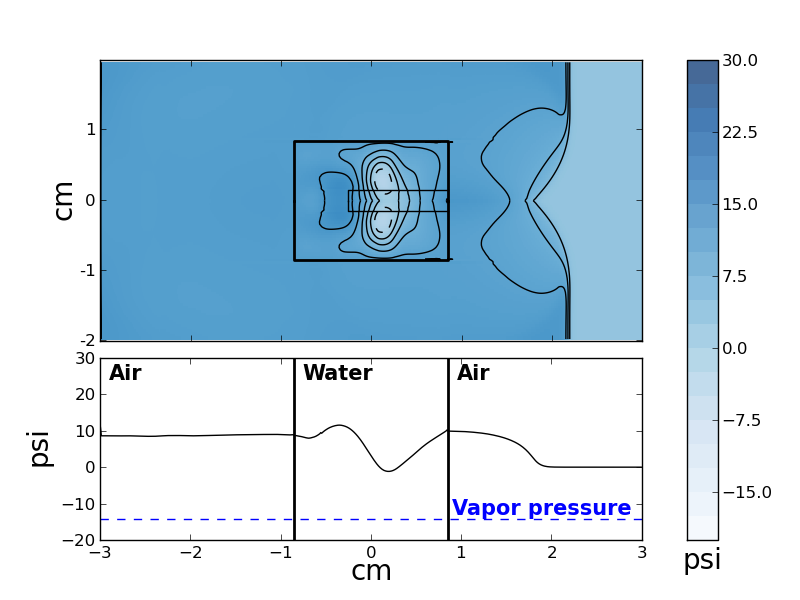} 
          \label{fig:TBIanim2D2_hyd}
  \caption{Same as \Fig{TBIanim3D} but with an hydrophone inserted. 
  Note that in Frame 4 the pressure does not go below the vapor pressure in this case.
  }
  \label{fig:TBIanim3D_hyd}
\end{figure}

The results from the simulation are shown for different times in \Fig{TBIanim3D} as contour and
pseudo-color plots of the pressure in the two-dimensional cross section. 
The corresponding one-dimensional pressure profiles along the axis of rotation are shown in the lower figure of each frame. Several 
relevant effects can be observed. The amplitude of the pressure is increased as expected from the 
previous one-dimensional calculations \cite{delrazo2014}. Also, we can see
that the geometry affects
the pressure profile as well as the ongoing reflections inside the cylindrical transwell. Of particular 
interest is the fourth time frame of \Fig{TBIanim3D}, where the reflected wave has a pressure below water
vapor pressure at room temperature. Since the water at room temperature can become gas when the pressure is
below the vapor pressure, cavitation is possible. It is known that cavitating bubbles can be 
responsible for cell detachment and cell membrane poration \cite{ohl2006sonoporation,chen2011observations}
and could be a possible mechanism of injury to the endothelial cells of the BBB.  

\begin{figure}[t]
  \centering
  \subfigure[]{ 
        \includegraphics[width=0.35\textwidth]{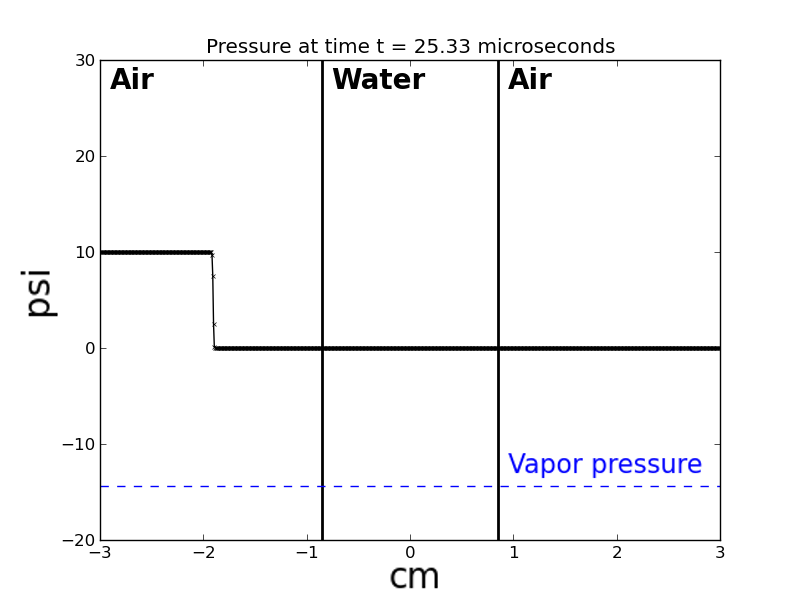} 
            \includegraphics[width=0.35\textwidth]{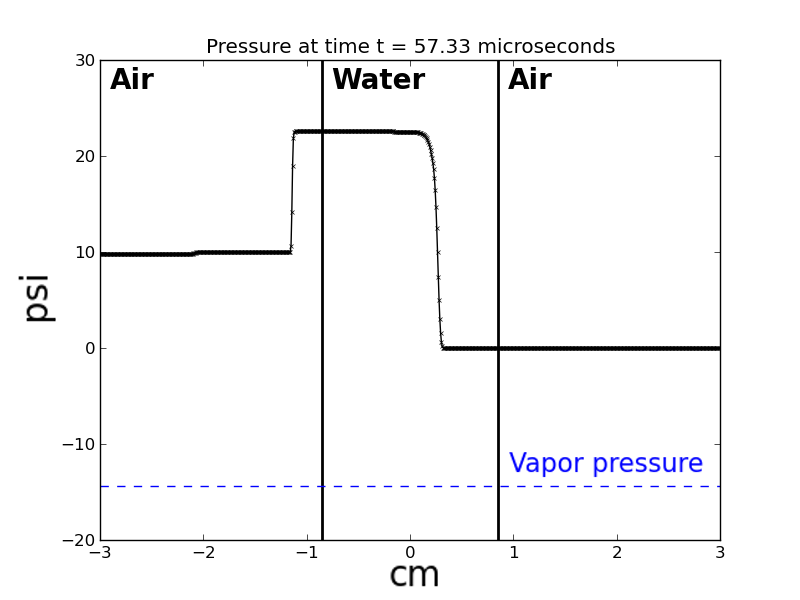} 
              \label{fig:TBIanim1D} }
 \subfigure[]{
        \includegraphics[width=0.35\textwidth]{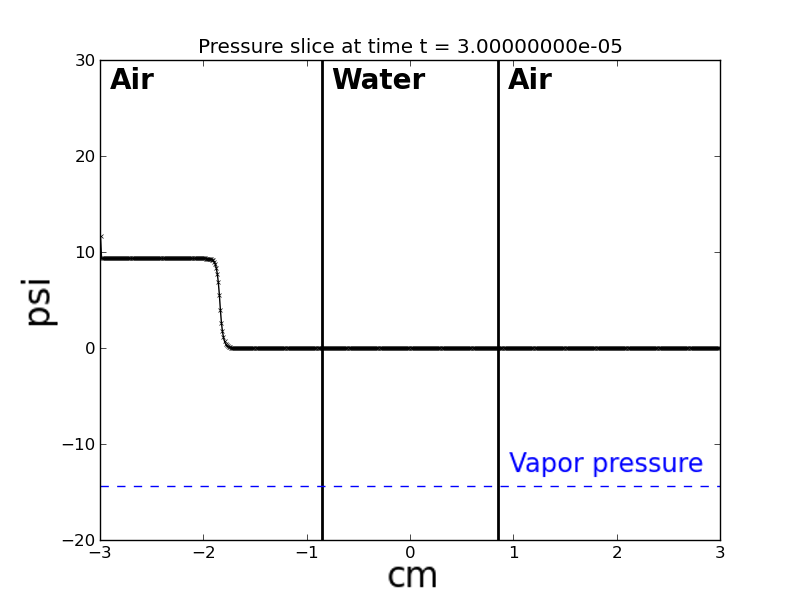} 
            \includegraphics[width=0.35\textwidth]{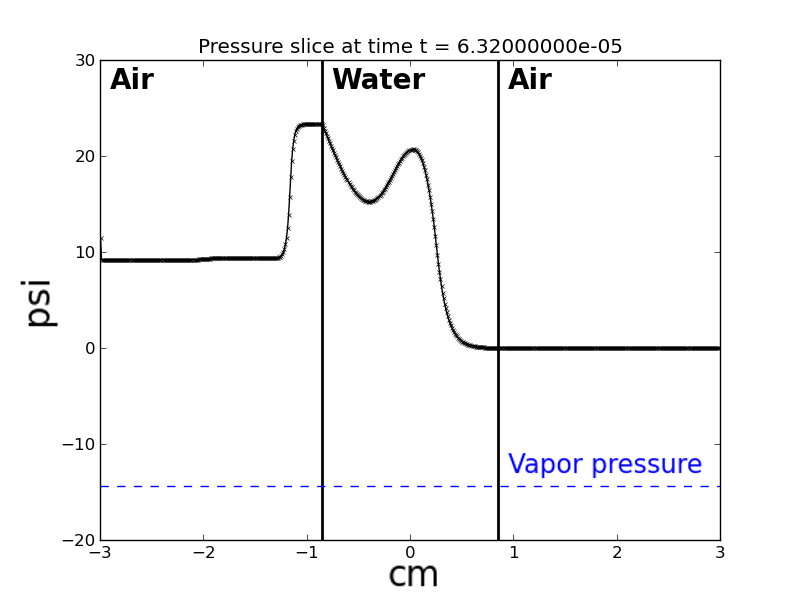} 
          \label{fig:TBIanim3Dslice} }
  \caption{
  (a) Pressure shown at two time frames from a one-dimensional simulation.
  Left:  The initial shock approaching the interface.
  Right: The reflected and transmitted shocks.
  (b) Pressure along the axis at the same two times, from the two-dimensional axisymmetric simulation.
  The edge effects in the pressure profile are evident in the second time frame. 
}
  \label{fig:TBIanim3Dn1D}
\includegraphics[width=0.88\textwidth]{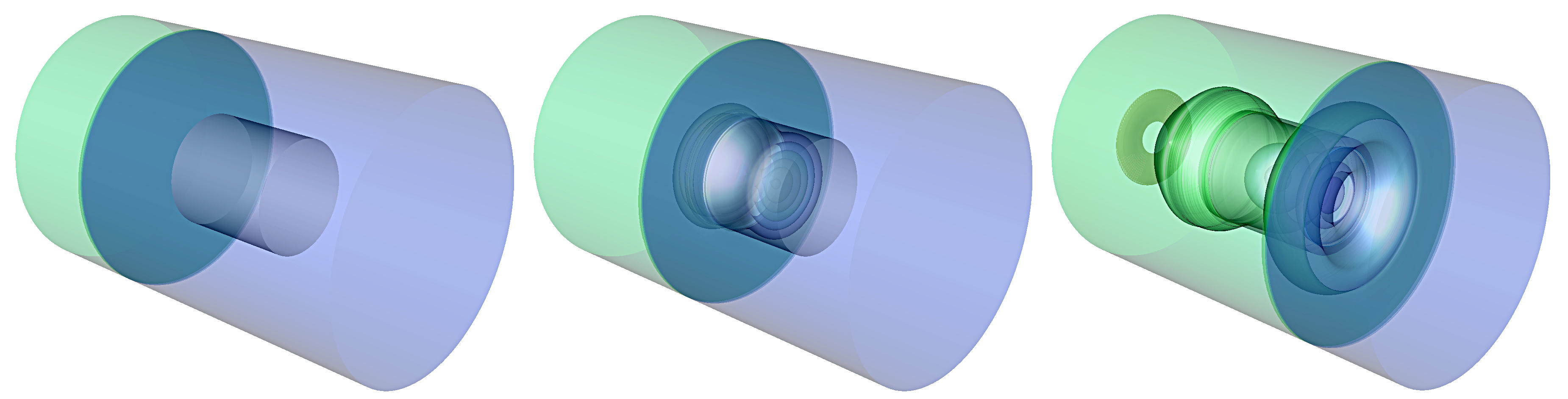}
\caption{Three-dimensional visualization by revolving the solution of frames 1,  3 and 6 of the two-dimensional axisymmetric
results from \Fig{TBIanim3D}. The cylindrical transwell can be well appreciated on the first frame. The visualization
shows the pressure contours, darker contours correspond to higher pressure. Its purpose is to emphasize that the two-dimensional
axisymmetric model is effectively modeling three-dimensional wave propagation.
}
\label{fig:COMPstube}
\end{figure}

To further understand these effects, we can observe 
\Fig{TBIanim3Dn1D} where the  axisymmetric
model is compared to the one-dimensional one. The geometrical edge 
effects are clearly seen in the second frame, where the pressure profile exhibits a decay
in the amplitude after the shock wave has crossed the interface.
This is  due to the presence of the cylindrical transwell walls parallel to
the axis of rotation.  As noted elsewhere \cite{delrazo2014}, pressure values below 
atmospheric pressures do not appear in the one-dimensional case, illustrating that low pressure 
values that might produce cavitation are a direct consequence of the geometrical edge effects.

As we mentioned before, we are employing a two-dimensional axisymmetric computational model, which effectively models
three-dimensional shock wave propagation. In \Fig{COMPstube}, we show a three-dimensional visualization of the
solution by revolving the solution of frames 1,3 and 6 of \Fig{TBIanim3D}. The figure shows three-dimensional
pressure contours, and it is included to emphasize the fact that we are modeling propagation of waves in three dimensions.

\subsection{Effects of introducing a hydrophone}
One might like to experimentally measure the pressure at the location of the endothelial cells in the
transwell in order to determine the force applied to the membrane and the possibility of cavitation.  
We attempted to introduce a customized version of the Y-104 hydrophone 
(Sonic Concepts, Bothell WA) 
in some of our laboratory experiments, but we were unable to gather 
sufficiently high quality low-frequency data to compare with our numerical results.
We did not pursue these experiments because we realized that the
introduction of this device could directly affect the signal being measured,
reducing the value of such data.  
A significant advantage of the computational model is that we can measure the pressure at 
computational gauge locations without interfering with the wave propagation.

\begin{figure}[t]
  \centering
  \subfigure[Gauge positions]{ \includegraphics[width=0.4\textwidth]{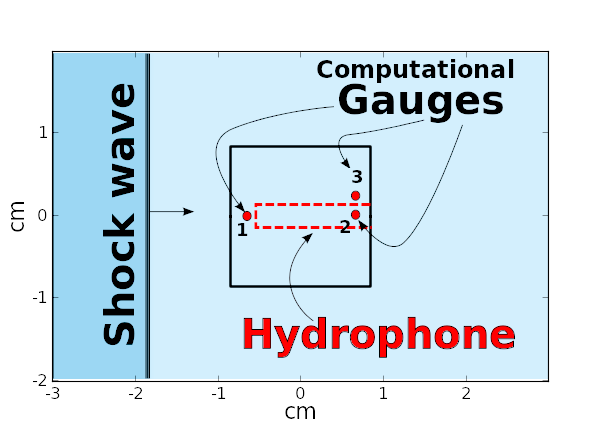}\label{fig:TBIdiag} }
  \subfigure[Gauge 1]{ \includegraphics[width=0.4\textwidth]{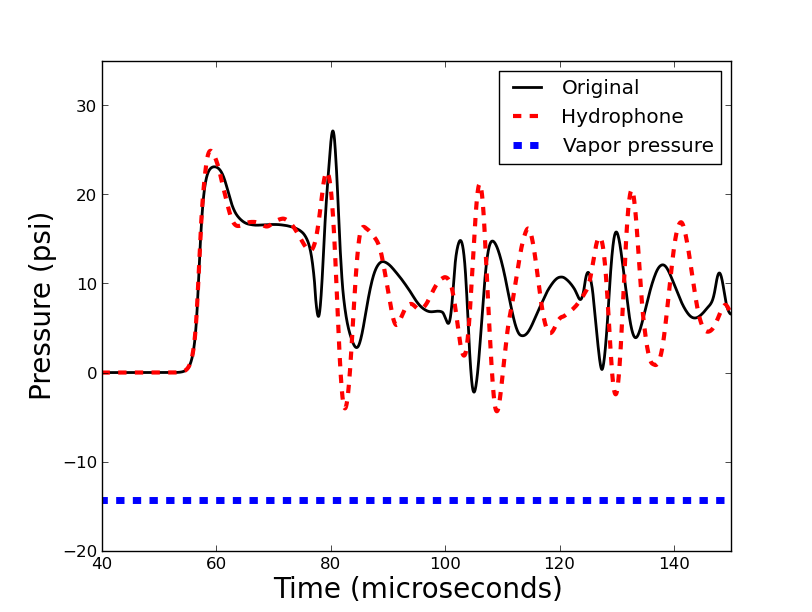} 
  \label{fig:TBIgauge1} }
  \subfigure[Gauge 2]{ \includegraphics[width=0.4\textwidth]{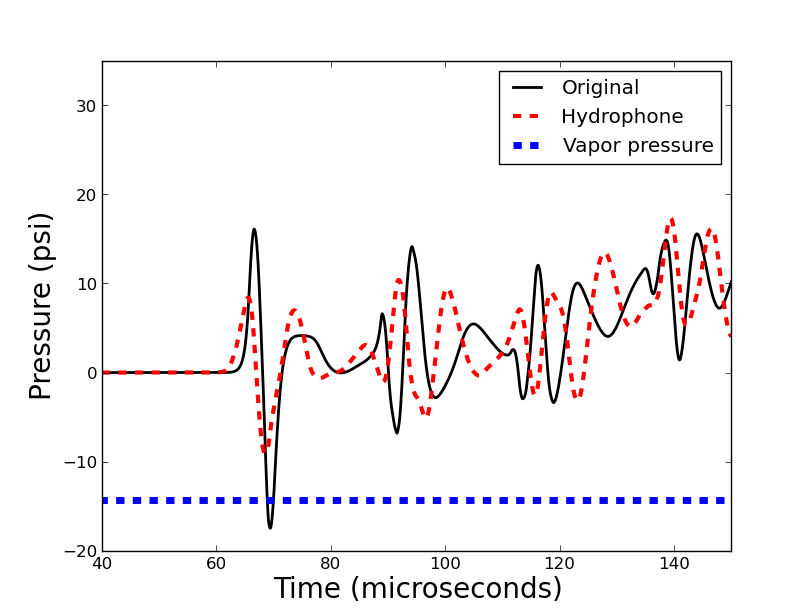}
   \label{fig:TBIgauge2} }
  \subfigure[Gauge 3]{\includegraphics[width=0.4\textwidth]{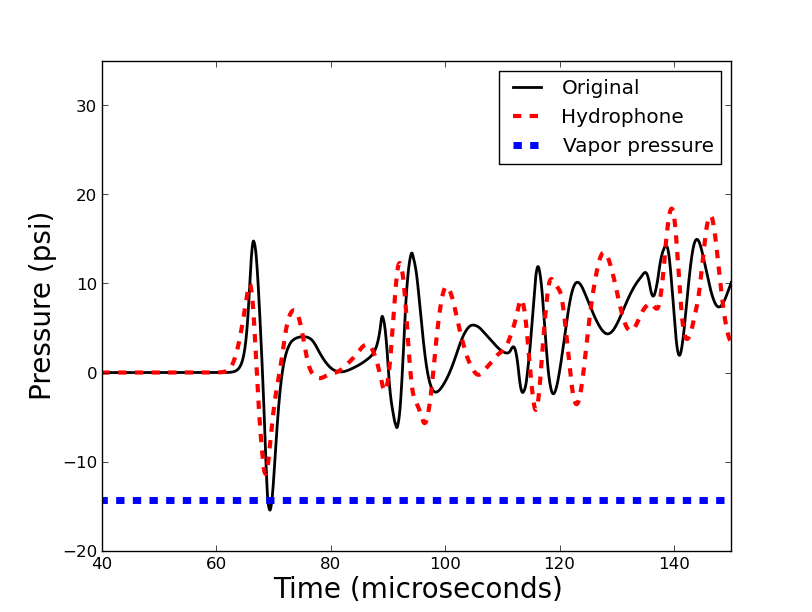} \label{fig:TBIgauge3} }
  \caption{Comparison of the pressure at computational gauges when a hydrophone is introduced with
  the pressure in the absence of a hydrophone. 
  The location of three gauges is shown on the first frame. The pressure profiles (psi) as a function of time 
  ($\mu$s) are shown for the three gauges. The output of the original simulation without the
  hydrophone is plotted with a solid line; the output of the simulation with the hydrophone 
  is plotted on a dashed line and the vapor pressure is plotted with a thick dashed 
  line. Note that the pressure falls below vapor pressure in the original simulation at Gauge 2 and Gauge 3 
  but not when the hydrophone is introduced. Also note that in the presence of the hydrophone, Gauge 2
  becomes irrelevant.}
  \label{fig:TBIgauges}
\end{figure} 

We can use the computational model to gain insight on how much the introduction of an hydrophone 
would change the experimental results.  To this end, we
include an axisymmetric computational hydrophone down the center of the
transwell in the following simulations, with a diameter of 2.85mm to match
the Y-104 model. The main effect that concerns us when incorporating the hydrophone in the
simulation is the reflection of acoustic
waves back into the liquid. The hydrophone is not uniquely 
composed of a single material and it is designed to have a net impedance of the same order 
of magnitude as water ($\sim 1.5 \times 10^6 Pa\cdot s/m$). Furthermore, solids usually have
impedances higher than water, so we can simply model the hydrophone as a general elastic solid with 
such properties. For this work, we model it as made of polystyrene with the 
parameters from Table \ref{tab:param} and a resulting impedance of ($\sim 2.4 \times 10^6 Pa\cdot s/m$). 
Modifying the impedance of the hydrophone material in the simulations through nine values within the same order of 
magnitude ($2 \times 10^6 Pa\cdot s/m$ to $7 \times 10^6 Pa\cdot s/m$) did not change any of the 
qualitative results presented here.

The computational results with the hydrophone are shown in \Fig{TBIanim3D_hyd}. We note there is a significant
difference between the results obtained in comparison to those without the hydrophone from \Fig{TBIanim3D}.
These data indicate that, in principle, hydrophone and intracranial pressure sensors placed in a small enclosed 
volume can alter shock wave propagation in functionally significant ways.
This has implications also for rodent experiments, as we see that an
intracranial pressure sensor placed within a volume comparable to that of a 
rodent skull can significantly alter shock wave dynamics, sufficient to change conditions that may favor cavitation. 

In order to better quantify the difference between the experiment with and without the hydrophone, we placed 
three gauges at key points in both systems. In \Fig{TBIgauges}, we can observe the comparison between the pressure 
profile as a function of time in the three chosen points. We can see the pressure only falls below 
vapor pressure in Gauge 2 and Gauge 3 when the hydrophone is not present. We can conclude that the inclusion of an 
hydrophone in the experimental system eliminated the possibility of observing cavitation. More importantly, 
measuring the pressure profile with a hydrophone in an experimental system like this one, affects the 
observed pressure profile, which supports the use of a computational model for quantifiable insight 
and answers to some experimental issues.

\section{Mathematical and computational models}
\label{sec:nummethods}

In this section, we give an outline of the numerical implementation,
summarizing the general methods used in Clawpack as well as the
original approaches  and implementations that were designed
uniquely for this work.

\subsection{The Euler equations}
\label{sec:TBInumimp}

We use the inviscid Euler equations for compressible flow, with different 
parameters in the equations of state (EOS) for each material.
The axisymmetric Euler equations in cylindrical coordinates ($r$,$\theta$,$z$) take the form
\begin{gather}
\label{eq:Eulercyl}
\begin{gathered}
  \fr{\p}{\p t}
    \left[\begin{array}{c} \rho \\   \rho u_r \\    \rho u_z \\    E   \end{array} \right] 
+ \fr{\p}{\p r} 
    \left[\begin{array}{c} \rho u_r \\   \rho u_r^2 + p \\    \rho u_r u_z \\     u_r(E+p)   \end{array} \right] + 
 \fr{\p}{\p z}
    \left[\begin{array}{c} \rho u_z \\   \rho u_r u_z \\    \rho u_z^2 + p \\     u_z(E+p)   \end{array} \right] 
=   \left[\begin{array}{c} -(\rho u_r)/r \\   -(\rho u_r^2)/r  \\    -(\rho u_r u_z)/r \\     -u_r(E+p)/r   \end{array} \right] ,
\end{gathered}
\end{gather}
where $\rho$ is the density; $u_r$ and $u_z$ denote the velocities in the radial and axial direction, $r$ and $z$ respectively; 
$E$ is the total energy and $p$ is the pressure. These equations have the same form as the two-dimensional Euler equations
with the addition of geometrical source terms (right hand side).  These source terms are further
discussed in Section~\ref{sec:srcterms}.

For the computational model, we must handle wave propagation
in liquid and elastic solids as well as in the air.  To handle this
range of materials we use the stiffened gas equation of state
(SGEOS), also known as the Tammann EOS. This equation of state is
very useful to model a wide range of fluids even in the presence
of strong shock waves and was successfully used in \cite{Fagnan2012, Fagnan:phd}
to model shock wave propagation in tissue and bone.
The Tammann EOS is given by
\begin{align}
  p = (\gamma - 1) \rho e -\gamma p_{\infty},
  \label{eq:SGEOS}
\end{align}
where $\gamma$ and $p_{\infty}$ can be determined experimentally for different materials and conditions. The choice of 
parameters for some materials is shown in Table \ref{tab:param}.
It is worth mentioning that for sufficiently weak shocks the 
Tammann EOS can be further simplified to the Tait EOS, which neglects
the energy coupling. In \cite{Fagnan2012} this was shown to be
adequate for modeling shocks in fluids and solids in the context of shock wave therapy. 
In this work, we will employ the Tammann EOS since it provides a more comprehensive approach
and conserves the energy coupling that could be useful to relate to thermodynamic quantities.

\begin{table}[H]
  \centering
  \begin{tabular}{| l || c | c |}
    \hline
     Material                & $\gamma$     & $p_\infty (GPa)$ \\ \hline 
     Air (Ideal gas EOS)     & 1.4          & 0.0             \\ \hline
     Polystyrene             & 1.1          & 4.79            \\ \hline
     Water                   & 7.15         & 0.3             \\
    \hline
  \end{tabular}
  \caption{Parameters for the Tammann EOS to model the different materials. The parameters for air and water were taken 
  from \cite{Fagnan2012}. Since the polystyrene is a solid, $\gamma$ was chosen very close to 1, and $p_\infty$ 
  was adjusted to yield the right speed of sound in polystyrene. The saline
solution in the transwell should have parameters very close to water.}
  \label{tab:param}
\end{table}

\subsection{Numerical methods}
\label{sec:FVMmethods}

The Euler equations (\ref{eq:Eulercyl}) are a hyperbolic system of conservation laws, so they can
be solved employing finite volume methods (FVM). This is done by using the wave propagation 
algorithms described in detail elsewhere
\cite{randysrbook,leveque1997wave} and implemented in Clawpack \cite{clawpack}. 
The fundamental problem that needs to be solved at each cell interface of our computation is the well 
known Riemann problem. A general one-dimensional Riemann problem for a system of conservation laws 
like Euler equations can be stated as
\begin{align}
  q_t& + f(q)_x = 0, \label{eq:rphll}\\
  q&(x,0) = \left\{ 
  \begin{array}{l l}
   q_\ell & \  \text{if $x<0$} \\
   q_r & \  \text{if $x>0$}, 
  \end{array} \right. \nonumber
\end{align}
where $q$ is the vector of conserved variables, $f(q)$ the corresponding 
fluxes and $q_\ell$ and $q_r$ constant states. 

When employing finite volume methods, we need to introduce the concept of cell 
average: $Q_i^n=\fr{1}{\Delta x}\int_{x_{i-1/2}}^{x_{i+1/2}} q(x,t_n)dx$,
where $i$ is the cell number and $n$ the time step index. At the edge between two cells, 
the Riemann problem initial condition would be determined by $q_\ell = Q_{i-1}^n$ and $q_r = Q_{i}^n$. 
After solving the Riemann problems at every cell edge, we can average the respective contributions 
to obtain the new cells average after a time $\Delta t$. The reader is referred elsewhere 
\cite{leveque1997wave,randysrbook} for a detailed exposition of the algorithms. 

The equations of motion are solved by 
implementing a hybrid Riemann HLLC-exact Riemann solver for the Euler equations 
with interfaces. This solver couples a Eulerian HLLC (Harten-Lax-van Leer-Contact) 
approximate Riemann solver, see \cite{toro1994restoration,torosbook} to a Lagrangian 
exact Riemann solver for the Euler equations with a Tammann EOS \footnote{A Lagrangian version of 
the HLLC solver can be also used}. As the interfaces are represented by contact discontinuities, the 
HLLC solver is ideal to deal accurately with interface problems. The method can be extended to two and 
three dimensions, retaining second order accuracy, by implementing transverse solvers with an unsplit 
method \cite{leveque1997wave}. We 
designed the transverse Riemann solvers as approximate solvers based on linear acoustics and adapted them to 
deal with interfaces. The source terms for the axisymmetric case are resolved using an operator 
splitting \cite{randysrbook,randysbbook}. A detailed description of the hybrid HLLC-exact normal Riemann 
solver for the Euler equations with the Tammann EOS with discontinuous
parameters is presented in 
\cite{delrazo2014}, in the context of one-dimensional problems. The extension
of this solver to a Riemann solver normal to a cell interface in two space
dimensions is straightforward and will not be discussed here.
For the unsplit wave propagation algorithms implemented in Clawpack, this
must be augmented with a transverse Riemann solver, as described in the Section \ref{sec:transverseRS}.  
The source terms that arise from axisymmetry are handled via a
fractional step approach described in Section~\ref{sec:srcterms}.

\subsubsection{Verification}
In this section, we will verify that the finite volume methods coupled with the hybrid Riemann HLLC-exact 
Riemann solver for the Euler equations with a Tammann EOS give the correct solution for a simple model problem. 
As the studies in Section \ref{sec:compresults}
are concerned with the dynamics of a shock wave traveling in an air-water-air system with two interfaces, we use this 
example as the test case. The exact analytic solutions of Riemann problems for Euler equations are only available
in one dimension, so we restrict our verification to a one-dimensional test.
This analysis will test the accuracy of the approximate hybrid Riemann solver, the key 
ingredient of our numerical method, also in the two-dimensional extension of the algorithm.

\begin{figure}
\centering
\includegraphics[width=0.6\paperwidth]{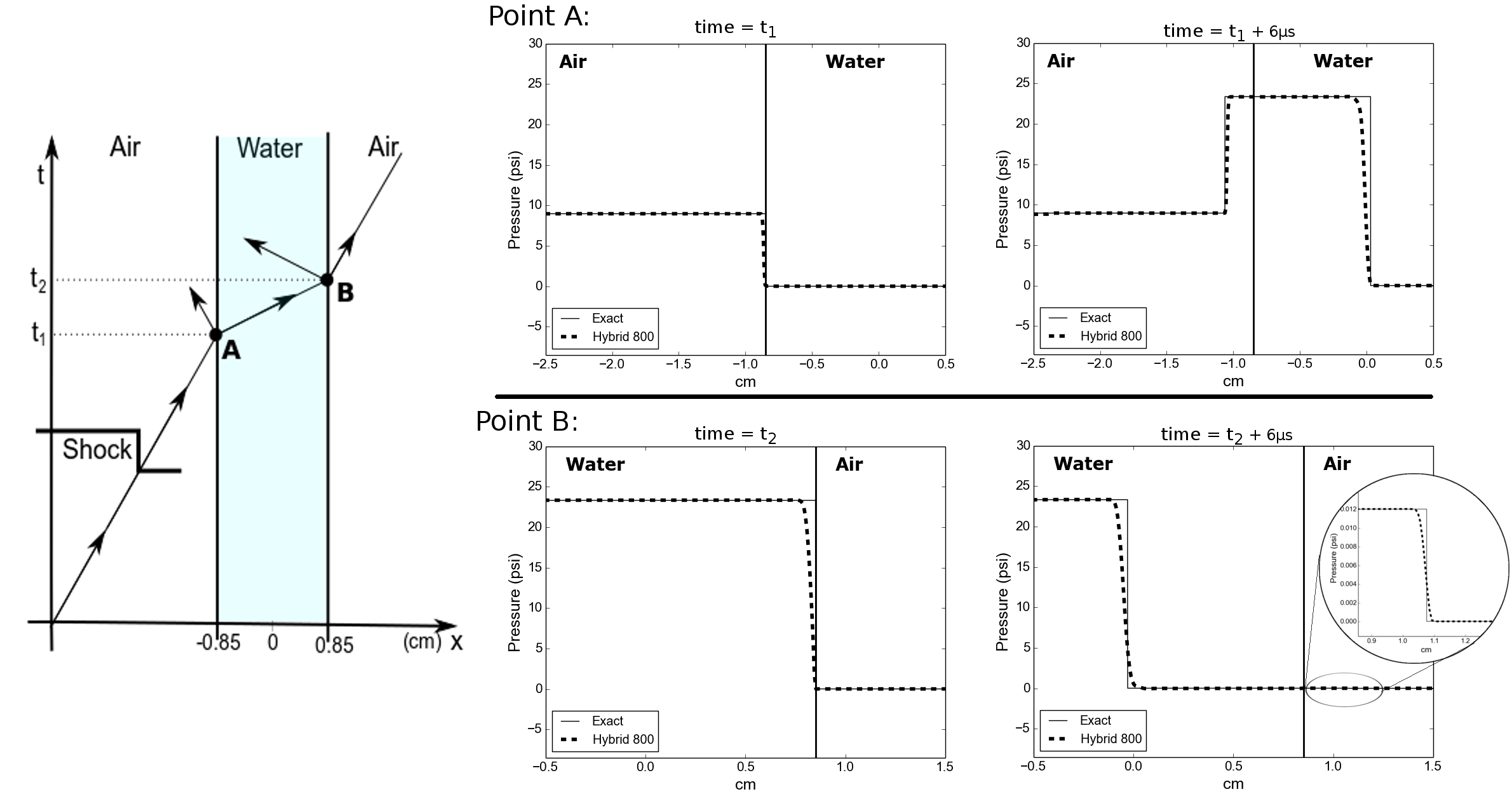}
\caption{Verification study of the numerical methods. On the left, we show the $x$--$t$ plane of an incoming shock
that hits the air-water interface at point A.  The transmitted shock then hits the water-air interface at point B. 
For both points, we show the solution when the shock hits the interface
and another one 6$\mu s$ later.
The exact solution can be constructed as described in the text and is compared with the numerical solution
computed using 
the hybrid HLLC-exact solver with 800 grid points. The air and water materials were modeled using the parameters 
from Table \ref{tab:param}.}
\label{fig:validation}
\end{figure}

The test problem is illustrated in the $x$--$t$ plane diagram on the left of \Fig{validation}. 
We will use a one-dimensional version of our algorithms, where we have an incoming shock of the same 
shape and intensity than the one used for \Fig{TBIanim3D}. We divide the domain into three materials: 
air-water-air, as in the original problem. At the time the shock hits the air-water interface at point A,
we can view the problem as a Riemann problem and compare it to the exact solution. Furthermore,
after the transmitted shock travels to the second interface, we have a second Riemann problem and
can repeat the same procedure at the water-air 
interface at point B and also compare the transmitted and reflected waves to the exact solution at that point. 
However, it should be noted that in the numerical algorithm
the incoming shock is not perfectly sharp, so we cannot expect a perfect match between our numerical solution
and the exact solution.

At point A in \Fig{validation}, we provide two plots: one just before the shock hits the air-water interface at $t_1$ where we
can frame the problem as a Riemann problem, and the second one 6$\mu s$ later. Both plots show two curves, one 
using the exact Riemann solver for the Euler equations with the Tammann EOS, with a 
jump in the parameters \cite{delrazo2015_02,Ivings1998}, and the other one is the numerical solution using 
the hybrid HLLC-exact solver for the Riemann problems that arise at each cell interface every time
step.

\begin{figure}
\centering
\includegraphics[width=0.45\paperwidth]{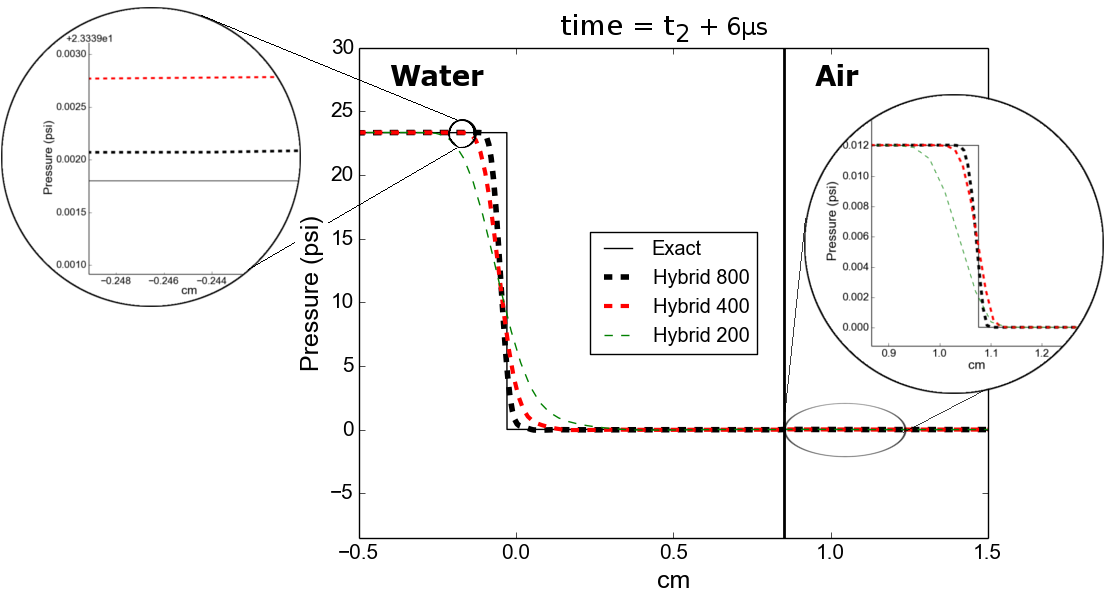}
\caption{The last plot in \Fig{validation} is recomputed with three different numerical resolutions
to show convergence. The numerical solutions are computed using the hybrid HLLC-exact solver with 200, 400 and 800 grid 
points. Two zoomed in regions in key areas are shown.}
\label{fig:convergence}
\end{figure}

The same procedure is repeated in point B of \Fig{validation}. The first plot shows the transmitted shock
from the exact Riemann solution at point A,
just before hitting the water-air interface at time $t_2$, and the second one 6$\mu s$ later. Notice in this last plot 
that there appears to be no transmitted wave.  However, the zoomed-in bubble shows that there is a very weak
transmitted shock at this interface of magnitude roughly 0.013 psi. 
Due to the much lower density of air relative to water, this interface acts nearly like a free boundary and the
reflected wave is a rarefaction wave, which is difficult to appreciate from
the figure since the difference between the rarefaction head and tail speeds is very small.
At both points in \Fig{validation}, we can see a very good agreement between our numerical solution and the exact solution. 
Furthermore, in \Fig{convergence}, we provide a convergence test. We chose to show it using the second plot at point B 
of \Fig{validation} since it gathers information of transmitted and reflected waves from both interfaces. 
\Fig{convergence} shows that the solution converges as we refine the resolution of our numerical solution.

\subsection{Transverse Riemann solvers}
\label{sec:transverseRS}
In order to obtain second order accuracy and improve stability in two-dimensional hyperbolic problems, 
the notion of a transverse Riemann solver was 
introduced in \cite{leveque1997wave}. This solver takes the results of a
Riemann solution in the direction normal to a cell interface and splits it
into components
moving in the transverse direction that contribute to updating the
solution in the adjacent rows of grid cells.  
Other alternatives also exist for solving multi-dimensional
conservation laws that attempt to use more fully multi-dimensional Riemann
solutions, for example in work of Roe \cite{roe:multid} and Fey \cite{fey:mot1,fey:mot2}.
Of particular relevance to the approximate Riemann solver approach used here
is the work of Balsara  \cite{balsara2012two,balsara2014multidimensional}, 
who defines two-dimensional HLLC Riemann solvers
that accept four input states that come together at an edge and outputs the 
multi-dimensionally upwinded fluxes in both directions. A comparison between
these two approaches could be of relevance in future studies.

For the present problem with sharp interfaces between very different
materials, instabilities were seen to easily arise, particularly at the
corners of the rectangular region representing the transwell. 
A special transverse solver was developed that we now describe, based on the
solver for acoustics in a heterogeneous media that is described in Section
21.5 of \cite{randysrbook}. Note that for two-dimensional problems on
rectangular grids, the cell average is calculated as
$Q_{i,j}^n=\fr{1}{\Delta y\Delta x}\int_{C_{i,j}} q(x,y,t_n)dxdy$,
where $C_{i,j}$ is the cell $[x_{i-1/2}, x_{i+1/2}]\times [y_{j-1/2},y_{j+1/2}]$.

\begin{figure}[t]
\centering
\includegraphics[width=0.38\paperwidth]{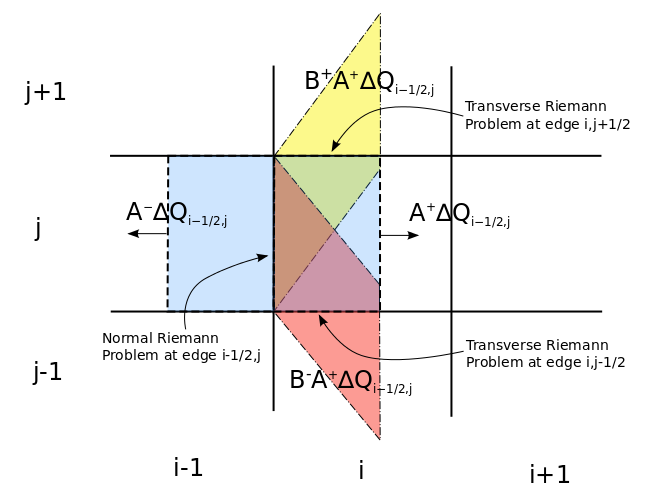}
\caption{Transverse solvers diagram for computational grid cells. The left-going and right going fluctuations of 
the normal Riemann problem at the edge between grid cells $(i-1,j)$ and $(i,j)$ is shown. The right-going 
fluctuation $A^+\Delta Q_{i-1/2,j}$ is decomposed into the up-going fluctuation $B^+A^+\Delta Q_{i-1/2,j}$
and the down-going fluctuation $B^-A^+\Delta Q_{i-1/2,j}$ by employing transverse Riemann solvers.} 
\label{fig:trans}
\end{figure}

We recall the basic idea of a transverse solver in \Fig{trans}.
For a constant coefficient linear hyperbolic system of equations 
$q_t + Aq_x +Bq_y =0$, the jump in normal flux between adjacent cells,
$A\Delta Q_{i-1/2} = A(Q_{i,j} - Q_{i-1,j})$, is split via the normal Riemann
solver into ``fluctuations'' $A^-\Delta Q_{i-1/2}$ and $A^+\Delta Q_{i-1/2}$
that correspond to the net contribution of all left-going or right-going
waves to the cell averages on either side.  Here $A^{\pm} = R\Lambda^{\pm}
R^{-1}$ where $A=R\Lambda R^{-1}$ is the eigen-decomposition of $A$ and
$\Lambda^\pm$ are the diagonal matrices in which either the negative or
positive eigenvalues have been set to zero.  Each fluctuation, e.g.
$A^+\Delta Q_{i-1/2}$, is then
further split into down-going and up-going components $B^-A^+\Delta
Q_{i-1/2}$ and $B^+A^+\Delta Q_{i-1/2}$, based on the matrices
$B^+$ and $B^-$.  

\def\A{{\cal A}}
\def\B{{\cal B}}

In the case of variable coefficients or nonlinear problems, 
the general notation $\B^-\A^+\Delta Q_{i-1/2}$ and $\B^+\A^+\Delta
Q_{i-1/2}$ is used for these two vectors.  For variable coefficient
acoustics, as described in \cite{randysrbook}, the up-going fluctuation from
the transverse splitting is based on eigenvectors of $B_{ij}$ and $B_{i,j+1}$,
while the down-going fluctuation is based on eigenvectors of $B_{ij}$ and
$B_{i,j-1}$.
For a nonlinear problem $q_t + f(q)_x + g(q)_y =0$, the eigen-decomposition
of some averaged Jacobian $g'(q)$ is generally used for the transverse
Riemann solver.

The present problem involves both nonlinearity and varying material
properties.  
Since we are modeling the almost incompressible liquid
in a Lagrangian frame of reference \cite{delrazo2014}, the transverse
Riemann problem will mostly be concerned with the two acoustic
waves. 
In order to derive the approximate transverse
solver, we will rely on linearized acoustic equations around $\rho_0,
u_0$ \cite{randysrbook} in terms of the density and momentum,

 \begin{gather}
\label{eq:linacoust}
\begin{gathered}
    \left[\begin{array}{c} \rho \\   \rho u \end{array} \right]_t +
    \underbrace{\left[\begin{array}{c c} 0 & 1 \\  c^2 & 0   \end{array} \right]}_{ \tilde{B}(Q)}
    \left[\begin{array}{c} \rho \\   \rho u \end{array} \right]_y
= 0,
\end{gathered}
\end{gather}
where we use $y$ as the space variable to emphasize this is solved in the transverse direction, $c$ is the sound speed 
and $\tilde{B}(Q)$ can be understood as a lower dimensional approximation of the transverse Jacobian 
$g'(Q_0)$ for the Euler equations. Note we assumed $u_0=0$, which is equivalent 
to assume we are in a Lagrangian frame of reference. The eigenvectors of the Jacobian of the system are given by $[1,\pm c]$ 
and the eigenvalues by $\pm c$; however, when solving the transverse Riemann problem, we might have different materials and 
sound speeds in the cell above or below. Instead of evaluating the whole Jacobian in one state, as in a Roe linear 
solver \cite{roe1981approximate}, we will evaluate the eigenvectors according to their location. These will be given by 
$v_U=[1,c_{U}]$ for the upward acoustic wave and $v_D=[1,-c_{D}]$ for the downward acoustic wave with eigenvalues $c_U$ and $-c_D$. 
Here $U$ and $D$ refer to cells $(i,j+1)$ and $(i,j)$ when computing $\B^+\A^+\Delta Q_{i-1/2,j}$ and 
to cells $(i,j)$ and $(i,j-1)$ when computing $\B^-\A^+\Delta Q_{i-1/2,j}$.
The matrix of eigenvectors $R$ and its inverse are given by,
 \begin{gather*}
\begin{gathered}
R=
    \left[\begin{array}{c c} 1 & 1 \\  c_U & -c_D   \end{array} \right],
\end{gathered}
 \ \ \ 
 R^{-1} = \frac{1}{c_U + c_D}\left[\begin{array}{c c} c_D & 1 \\  c_U & -1   \end{array} \right].
\end{gather*}
The up-going and down-going fluctuations for $\A^+\Delta Q_{i-1/2,j}$ are obtained by expanding the fluctuation 
in terms of these two eigenvectors or waves, $\A^+\Delta Q_{i-1/2,j} = \alpha_U v_U + \alpha_D v_D$,
so we need to solve $R\alpha=\A^+\Delta Q_{i-1/2,j}$, 
Note that the required fluctuation $\A^+\Delta Q_{i-1/2,j}$ for 
the Euler equations is a 4 dimensional vector with fluctuations in density, normal momentum, transverse momentum and internal
energy. As we are only interested in the acoustic waves, we will assume the fluctuations in normal momentum and energy are 
negligible, so we define the acoustic part of the fluctuation as the first and third entry of the 4 dimensional vector, 
i.e. $\A^+_{ac}\Delta Q_{i-1/2,j} = [\A^+_{\Delta Q 1},\A^+_{\Delta Q 3}]$. Solving the system for 
the vector $\alpha= R^{-1}\A^+_{ac}\Delta Q_{i-1/2,j}$, we obtain
\begin{align*}
 \alpha_U = \frac{1}{c_U + c_D}\left( c_D \A^+_{\Delta Q 1} + \A^+_{\Delta Q 3}  \right), \\
 \alpha_D = \frac{1}{c_U + c_D}\left( c_U \A^+_{\Delta Q 1} - \A^+_{\Delta Q 3}  \right).
\end{align*}
The up-going and down-going acoustic fluctuations are given by the velocity times the waves,
\begin{align*}
 \B^+_{ac}\A^+\Delta Q_{i-1/2,j} = c_U \alpha_U v_U, \\
 \B^-_{ac}\A^+\Delta Q_{i-1/2,j} = - c_D \alpha_D v_D.
\end{align*}
We require to solve two of these transverse solvers for the Euler equations as shown in the grid in \Fig{trans}. 
We will only consider the up-going fluctuation of the transverse solver at $(i,j+1/2)$ and the down-going fluctuaction of the
solver at $i,j-1/2$. This yields the full fluctuations as
 \begin{gather*}
\begin{gathered}
\B^+\A^+\Delta Q_{i-1/2,j} = 
 \frac{c_3 \left( c_2 \A^+_{\Delta Q 1} + \A^+_{\Delta Q 3}  \right) }{c_3 + c_2}
\left[\begin{array}{c}  1 \\ 0 \\  c_3 \\ 0   \end{array} \right], \\
\B^-\A^+\Delta Q_{i-1/2,j} = 
 \frac{-c_1\left( c_2 \A^+_{\Delta Q 1} - \A^+_{\Delta Q 3}  \right)}{c_1 + c_2}
\left[\begin{array}{c}  1 \\ 0 \\  -c_1 \\ 0   \end{array} \right],
\end{gathered}
\end{gather*}
where $c_1,c_2$ and $c_3$ are the speeds of sound in cells $(i,j-1)$, $(i,j)$ and $(i,j+1)$ respectively and the 
non-acoustic fluctuations were neglected. The sound speeds are calculated with the pressure, density and the 
parameters of the Tammann EOS in the respective cell with $c=\sqrt{\gamma\frac{p+p_{\infty}}{\rho}}$. Note this process
is repeated in exactly the same manner for the left going fluctuation $\A^-\Delta Q_{i-1/2,j}$ of the normal
Riemann problem. 

\subsection{Geometrical source terms}
\label{sec:srcterms}
In order to solve for the source terms of equation (\ref{eq:Eulercyl}), we need to apply a splitting method, see 
\cite{randysrbook}. In the first half time step, we solve the homogeneous version of equation (\ref{eq:Eulercyl}) 
over the whole grid, and in the second step we solve the system of ODEs obtained by ignoring the flux terms,
\begin{gather}
\label{splitting}
\begin{gathered}
  \fr{d}{d t}
    \left[\begin{array}{c} \rho \\   \rho u_r \\    \rho u_z \\    E   
          \end{array} \right] 
=   \left[\begin{array}{c} -(\rho u_r)/r \\   -(\rho u_r^2)/r  \\    -(\rho u_r u_z)/r \\ -u_r(E+p)/r 
          \end{array} \right].
\end{gathered}
\end{gather}
This equation can be solved with any explicit time integrator method like forward Euler and Runge-Kutta methods
or an implicit solver, such as TR-BDF2. However, this particular system can be solved exactly. Consider the first equation
of equations (\ref{splitting}) and multiply it by $u_r$, then
\begin{align*}
u_r \fr{d\rho}{dt} &= \fr{\rho u_r^2}{r}, \\
\Rightarrow \fr{d\rho u_r}{dt} - \fr{d u_r}{dt}\rho &=  -(\rho u_r^2)/r,
\end{align*}
where we used the product rule. Now substituting the second equation of (\ref{splitting}) into this result, 
we obtain $\fr{d u_r}{dt}=0$, so $u_r$ is constant. The same procedure can be applied to obtain that $u_z$ is
also constant.

As the total energy is given by $E = \rho e + \frac 1 2  \rho(u_r^2 + u_z^2)$, where the Tammann EOS (\ref{eq:SGEOS})
allows the substitution $\rho e = (p + \gamma p_\infty) / (\gamma - 1)$. As $u_r$ and $u_z$ are constant, we can 
differente the energy, $E_t = (\rho e)_t = \frac{1}{\gamma -1} p_t$. These results in conjuction with the fourth 
equation of (\ref{splitting}), yield $p_t = -(u_r / r) [\gamma (p + p_\infty) + \frac 1 2  (\gamma-1)\rho(u_r^2 + u_z^2)]$.
We now have a full system of equations in the primitive variables:
\begin{align*} 
\fr{d \rho}{dt} &= -(u_r/r) \rho,\ \ \ \
\fr{d u_r}{dt} = 0,\ \ \ \
\fr{d u_z}{dt} = 0,\\
\fr{d p}{dt} &= -(u_r / r) \left(\gamma (p + p_\infty) + \frac 1 2
(\gamma-1)\rho(u_r^2 + u_z^2)\right).
\end{align*} 

The first three equations can easily be solved, and the fourth equation can also be solved
with the solution of the first one and an integrating factor. Using the fact that the initial
conditions for the computation are the variables at time $t^n$, and we want the solution
at time $t^{n+1} = t^n + \Delta t$, we obtain

\begin{equation} 
\begin{split} 
\rho^{n+1} &= \exp\left(-\fr{\Delta t u_r^n}{r}\right) \rho^n,\ \ \ \
u_r^{n+1} = u_r^n,\ \ \ \
u_z^{n+1} = u_z^n,\\
p^{n+1} &= \exp\left(-\fr{\Delta t \gamma u_r^n}{r}\right) p^n 
- p_\infty \left(1 - \exp\left(-\fr{\Delta t \gamma u_r^n}{r}\right)\right)
\\& \ \ \ -\frac {\rho^n}{2}\left((u_r^n)^2 + (u_z^n)^2 \right)
\left[\exp\left(-\fr{\Delta t u_r^n}{r}\right) - \exp\left(-\fr{\Delta t \gamma u_r^n }{r}\right)\right]
, \\
E^{n+1} &=  \fr{p^{n+1} + \gamma p_\infty} {\gamma -1} + \frac 1 2 \rho^{n+1}
\left((u_r^n)^2 + (u_z^n)^2\right).
\end{split} 
\end{equation}

The parameters $\gamma$ and $p_\infty$ are given by the Tammann EOS in equation (\ref{eq:SGEOS}). The 
equations we just obtained allow us to calculate one-time step of (\ref{eq:Eulercyl}) in our splitting method. Note 
these source terms are never singular in the computation; when using finite volume methods, the quantities are 
evaluated at cell centers, so $r>0$.

\section{Discussion}
\label{sec:disc}

A computational model was designed to better understand the physical
forces developed by blast-induced shock waves that can damage brain
endothelial cells in an {\em in vitro} model of the BBB. The numerical
modeling of the experiment employs finite volume methods and 
requires coupling a highly compressible material (air) with 
a nearly incompressible liquid contained
in a fixed region in space. The coupling is accomplished by employing
a Tammann EOS and designing both normal and transverse Riemann solvers
that can couple these two materials --- one in a Eulerian frame of
reference and the other in a Lagrangian frame of reference. Results
show the shock wave pressure amplitude and velocity increase when
crossing from air to the water (saline solution). This is in agreement with the one-dimensional
simulations described by us previously \cite{delrazo2014}, as well as other works
mentioned in a recent review \cite{gupta2013mathematical}. One aspect
of the potential relevance of this effect lies in the underestimation
of the pressure intensities experienced by the cells when one considers 
only the amplitude and kinetic properties of a standard open field blast overpressure.


Comparison of the computational results here to the one-dimensional tests
performed in \cite{delrazo2014} show that the transwell geometry is very
relevant. The edge effects from the cylinder, combined with the rarefaction wave
arising when the shock reflects off the distal end of the transwell,
can generate low enough pressure to potentially
produce cavitation, which could be a cause of cell damage
\cite{ohl2006sonoporation}.
The simulation with a hydrophone in place
does not show low enough pressure values to produce cavitating bubbles. 
These results indicate that the computational model could be useful
to experimentalists in analyzing how the
introduction of a measuring device affects the outcome of the experiment 
and the likelihood of cavitation being a BBB tissue damage mechanism.

Although based on an idealized
model, our computational approach allows us to measure the pressure profile at any point and at the exact location of the biological sample without
interfering with the actual experimental setup. This task would be extremely difficult to 
obtain empirically. The high-resolution signal obtained by our computational method allows 
us to apply it to identify regions with low enough pressure to potentially produce cavitation. 
Furthermore, our results allow us to suggest cavitation as a damage mechanism that might explain
the experimental results, for instance, the mislocalization of the tight junction proteins, ZO1, and claudin-5,
that functionally disturb the BBB. This kind of study can clarify the
qualitative behavior of the system and, where it is impossible for
experimentalists. It can also suggest possible connections between damage mechanisms 
and anatomical, functional, morphological, and molecular specificity, obtained from the experimental 
results.

The computational model developed in this work 
was designed for a specific application; however, the methods developed can be
adapted and applied to other experiments with similar simplified
geometry. These methods can also be extended to other geometries
and the Clawpack software (with adaptive mesh refinement)
can be applied in situations where a logically 
rectangular grid can be mapped to a quadrilateral
two-dimensional grid. 
This can include situations in which the interface is circular
or of other smooth shape lacking corners using the sort of mappings
proposed in \cite{cal-hel-rjl:circles}, which have been used for elastic and
poroelastic wave propagation problems in the work of Lemoine
\cite{lemoine:phd,LemoineOu2014}. Extension of the methods
proposed in this paper to such cases is currently under way 
and will be reported elsewhere \cite{delrazo2015_02}. This extension is 
clinically relevant; it allows detailed studies of the pressure signal 
obtained by shock waves interacting directly with the skull in conditions that 
might not be feasible experimentally,  emphasizing the importance 
of having a computational model available.

The computational simulations were evaluated up through the first
200 microseconds. As seen in \Fig{stube-data-c}, this corresponds
to a very short time period behind the shock, before the bulk of
the trailing rarefaction wave has passed the transwell.  Planned
future work includes the refinement of our numerical method to carry
out the simulation to longer times. This can be of relevance given
the negative pressure values and oscillations that arise on millisecond
time scales, as well as the secondary reflection-induced shock, see
\Fig{stube-data-c}. These features, along with the internal reflections
might also cause or even increase cavitation effects.

Some other possible future research directions include extension of the
computational methods to arbitrary interface geometry and to two-phase
models that can simulate cavitation. 
In addition, the {\em in vitro} system coupled with the
computational model can be used for future clinically relevant studies.
The ability to determine pressure traces at the precise location of the planar
endothelial cell monolayer could be used as an input into a mechanical
model of membrane dynamics during blast wave propagation. This would
permit new and highly refined estimates of the physical forces that
brain endothelial cells may be exposed to, such as high frequency 
BBB oscillations that may disrupt cellular functions even without 
gross brain displacements.

An important novel aspect
of this approach is that these estimates can be correlated to specific
quantifiable measurements of cellular damage, dysfunction of the BBB
as a system of interacting cells, and even aberrant subcellular
protein trafficking where it is possible to investigate the mechanisms
by which blast alters how critical BBB proteins, such as claudin-5
(Appendix, Figure 5.2) are misdirected inside cells away from tight
junctions. 

The simulation code developed in this work is available at
\cite{TBI-zenodo}, along with the raw data and SPSS statistical
analysis discussed below in Section~\ref{sec:permeability}.
The simulation code relies on Clawpack \cite{clawpack} and the results
presented in the paper were obtained with Version 5.2.2. 


\section{Appendix: Additional experimental results and methodology}
\label{app:methods}

Using well-established methods \cite{banks2004triglycerides,banks2005effects,nakaoke2005human}
mouse brain-derived endothelial cells (MBECs), purified from wild-type C57BL6 mice, were grown on 
permeable nylon support membranes in standard transwell chambers (see \Fig{stube-data-a}) and formed 
endothelial cell monolayer tight junctions that functionally mimic the BBB, which 
is responsible for maintaining and regulating separation between the central nervous system (CNS) and 
the circulating peripheral blood supply \cite{banks2010blood,zlokovic2011neurovascular}. The transwell 
chambers were filled completely with an aqueous solution (serum-free DMEM/F12 medium containing bFGF (1 ng/ml) 
and hydrocortisone (500 nM)). 
For blast exposure, the transwells were secured in the shock tube with the bottom of the transwell facing the
oncoming shock wave (see Figure 1.2). For all experiments, BBB cells were exposed to a single mild blast of indicated
intensity (psi).

\begin{figure}[t]
\centering
\includegraphics[width=0.38\paperwidth]{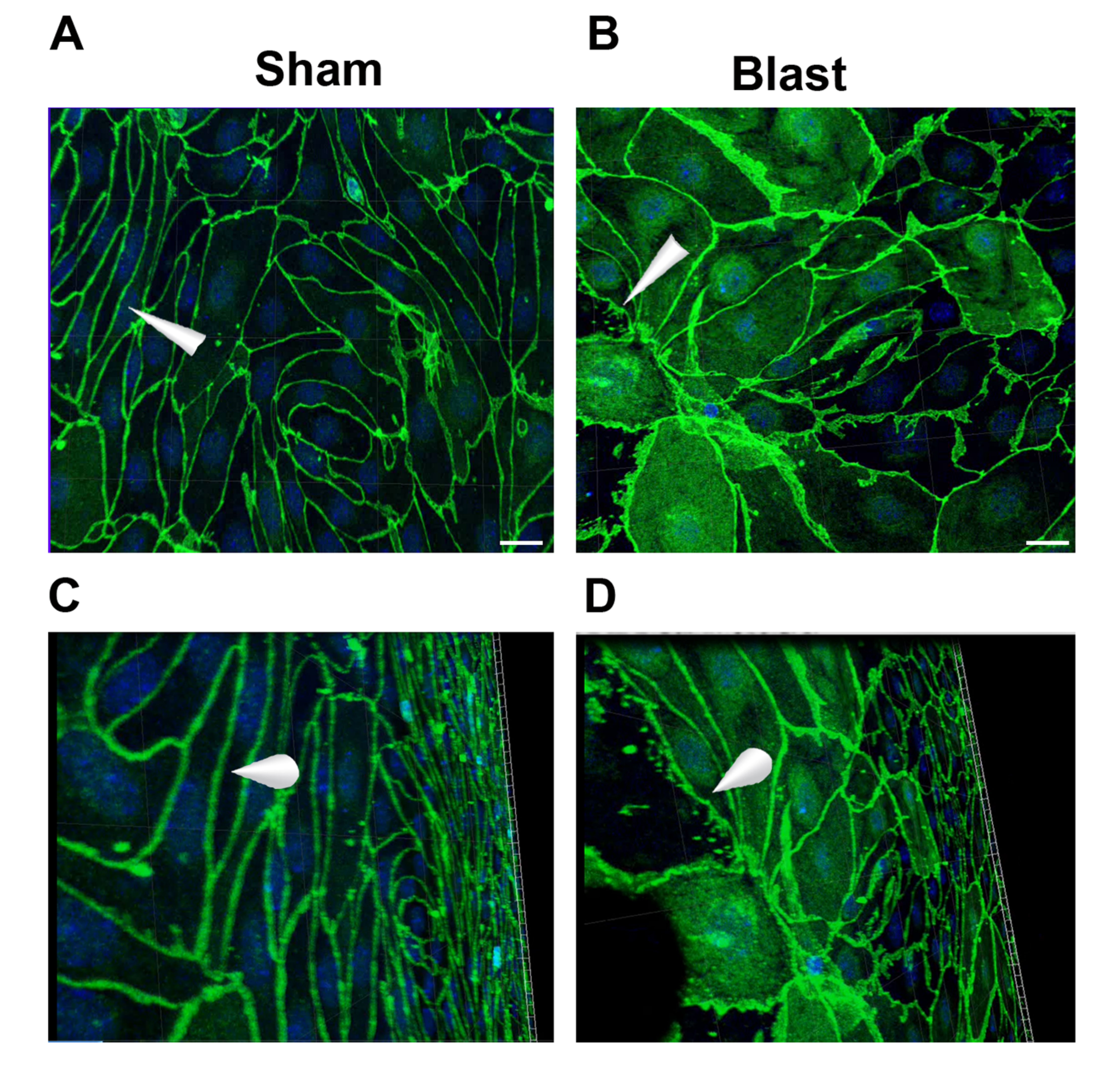}
\caption{
(A) Laser confocal microscopy reveals normal ZO-1 expression patterns expressed specifically at uniform, 
well-defined tight junctions along cell-to-cell interfaces within the plane of the brain-derived microvessel 
endothelial cell monolayer. (B) In contrast to the sham condition, ZO-1 expression in blast-exposed endothelial 
cells is highly dystrophic with widespread mislocalization in cellular domains remote from tight junctions. 
Panels A and B show a merged, serial reconstruction comprised of 27 images acquired at $0.2 \mu$m intervals 
along the z-axis orthogonal to the plane parallel with the MBEC cell monolayer. (C and D) Lower panels show 
oblique x-y-z plane views of the panels above (A,B), thereby permitting an improved assessment of blast-induced 
tight junction dysmorphology compared to normal sham tight junctions. Nuclei are stained blue with Dapi. 
Arrowheads denote the same cell-to-cell contact 
domains in the corresponding sham (A, C) and blast (B, D) images. Scale bars = $20 \mu$m.} 
\label{fig:ZO1}
\end{figure}

In addition to the experiment presented in Section \ref{sec:BBBexp}, we performed another experiment to investigate the 
effects of the shock tube blast exposure on tight junction morphology. Singly blasted (13-13.9 psi) or sham-treated 
monolayers were immunostained with antibodies recognizing the tight junction-associated scaffolding protein, ZO-1 
\cite{tokuda2014zo} 24 hours after treatment and then imaged using laser confocal microscopy. ZO-1 expression in 
sham-treated MBEC monolayers appeared morphologically normal with ZO-1 immunostaining tightly restricted to the 
interposing plasma membrane domains at points of cell-to-cell contact (\Fig{ZO1}A). In marked contrast
to this, blast exposure induced ragged, hypertrophic appearing tight junctions (\Fig{ZO1}B). In addition, ZO-1 expression 
appeared mislocalized in association with peri- abluminal and/or peri-luminal plasma membranes domains. This expression 
pattern is also consistent with diffuse intracellular cytoplasmic ZO-1 mislocalization.

The confocal images in \Fig{ZO1}A and \Fig{ZO1}B are maximum-field projections comprised of 27 merged images collected 
at $0.2 \mu$m step intervals in the z-axis orthogonal to the plane of the MBEC monolayer, thereby representing a total 
depth of $5.4 \mu$m that encompassed the full cross-sectional width of the MBEC monolayers. \Fig{ZO1}C and \Fig{ZO1}D 
depict three-dimensional serial reconstructions of images in the upper panels projected at oblique angles. For ease of 
reference, the arrowheads denote the same cell-to-cell contact points in panels A, C and B, D (sham and blast-exposed, 
respectively). From these oblique angles the degree of blast-induced tight junction dysmorphology and ZO-1 
mislocalization are more easily appreciated (Also, see supplementary videos). 

\begin{figure}[t]
\centering
\includegraphics[width=0.38\paperwidth]{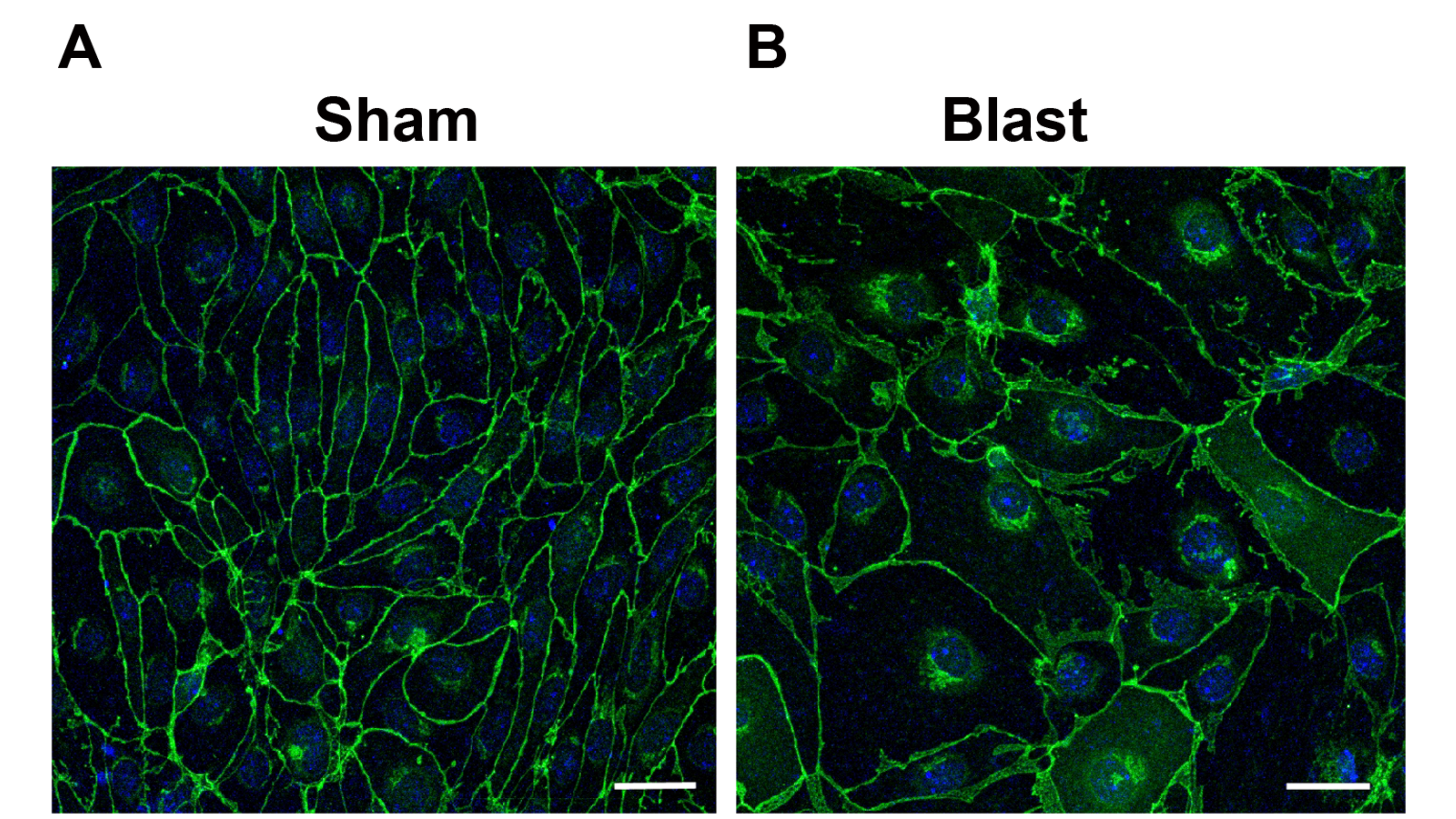}
\caption{
(A) Laser confocal microscopy reveals normal claudin-5 expression at tight junctions localized along cell-to-cell 
contacts of the MBEC monolayer. (B) In contrast to the sham controls, claudin-5 expression in blast-exposed 
endothelial cells is dysmorphic, indicative of aberrant tight junction structure.  In addition, claudin-5 is 
broadly mislocalized and accumulates in asymmetric peri-nuclear intracellular compartments, strongly suggesting 
that blast exposure induces aberrant subcellular trafficking of claudin-5. Nuclei are stained blue with Dapi. Scale bars = $25\mu$m.} 
\label{fig:Claudin}
\end{figure}

Claudin-5 is a tight junction-specific membrane bound protein \cite{jan2014molecular} that is a critical 
regulator of BBB permeability \cite{nitta2003size}. \Fig{Claudin} shows that a single mild blast exposure also markedly 
disrupted claudin-5 expression.  As with Z0-1, claudin-5 immunostaining revealed aberrant, hypertrophic appearing tight 
junctions in the blast-exposed monolayers. In addition, the asymmetric peri-nuclear claudin-5 immunostaining clearly 
demonstrates that blast exposure caused it to become aberrantly retained within the cells, thus raising the possibility 
that normal polarized subcellular trafficking of claudin-5 into and/or away from tight junction domains may be disrupted in the 
blast-exposed MBECs. 

These data suggest that blast exposure causes
mislocalization of the tight junction proteins,
ZO1 and claudin-5, away from tight junctions. Previous work using the continuous cell line, bEnd.3 showed that 
blast causes a loss of ZO1 and claudin-5 \cite{hue2013blood,hue2014repeated}. This difference could be because bEnd.3 cells are less 
differentiated than brain-derived microvessel endothelial cells, and which form barriers with lower TEER values 
than primary brain endothelial cultures used in this report \cite{deli2005permeability}. Nonetheless, our findings in BMECs and in vitro 
blast studies using bEnd.3 cells \cite{hue2013blood,hue2014repeated}, collectively demonstrate that blast exposure disturbs expression of 
proteins critical for maintaining BBB integrity.

Mechanistically, protein mislocalization suggests a dynamic 
alteration in the cellular process of adjusting to injury, whereas overall tight junction protein loss may suggest 
co-attending endothelial cell death or impaired protein production or increased tight junction proteolysis. 
Increasingly, tight junction protein mislocalization is viewed as an underlying pathology in diseases with BBB disruption 
and is the pattern, for example, in inflammatory conditions \cite{dimitrijevic2006effects,andras2003hiv}. 

\subsection{Culture of primary brain microvascular endothelial cells}
Brain microvascular endothelial cells (BMECs) were isolated from 6-8 week old CD-1 mice based on established standard
with some modifications procedures \cite{coisne2005mouse,jacob2010c5a}. All procedures involving animal subjects were 
carried out following protocols approved by the Veterans Affairs Puget Sound Health Care System Institution Animal Use 
and Care Committee (IACUC). Briefly, meninges were removed from freshly dissected brain cortices, and then the brain was minced. The minced 
brain matter was ground using a Dounce homogenizer in Dulbecco's Modified Eagle's Medium/Nutrient Mixture F-12 Ham 
(DMEM/F12; Sigma-Aldrich) supplemented with gentamicin ($50 \mu g/ml$; Sigma-Aldrich). 30\% Dextran (v/v; from 
Leuconostoc spp., MW ~70,000 Da; Sigma-Aldrich) was added to the homogenate 1:1 and supplemented with 10\% bovine serum 
albumin (BSA, Sigma-Aldrich) to achieve a final concentration of 0.1\%. The mixture was centrifuged at 3000 × g for 25 
$\mathrm{min}$ at $4\celsius$. The pellet obtained after the centrifugation was re-suspended in DMEM/F12, filtered through a
$70 \mu m$ nylon mesh, and centrifuged again at 1000 × g for 10 $\mathrm{min}$ at room temperature (RT). The resulting pellet
was digested at $37\celsius$ for 30 $\mathrm{min}$ with DMEM/F12 containing collagenase (0.2 U/ml), dispase (1.6 U/ml; 
collagenase/dispase, Roche Life Sciences) and DNase I ($10 \mu g/ml$; Sigma-Aldrich). The digested vessel suspension was 
filtered through a $ 21 \mu m$ nylon mesh. The filtrate was washed several times with DMEM/F12, and the resulting capillary 
suspension was seeded on dishes coated with collagen type IV (0.1 mg/ml; Sigma-Aldrich) and fibronectin (0.1 mg/ml; Sigma-Aldrich). 
BMECs were cultured in BMEC medium, consisting of DMEM/F12 supplemented with 20\% plasma-derived fetal bovine serum 
(Animal Technologies), 1\% GlutaMAX (Life Technologies), basic fibroblast growth factor (bFGF, 1 ng/ml; Roche Life Sciences), 
heparin ($100 \mu g/ml$), insulin ($5 \mu g/ml$), transferrin ($5 \mu g/ml$), selenium (5 ng/ml) (Insulin-transferrin-selenium 
medium supplement; Life Technologies), and gentamicin ($50 \mu g/ml$; Sigma-Aldrich. Puromycin (4 $\mu$g/ml; Sigma-Aldrich) was 
added to BMEC medium for the first 48 hours after plating to remove pericytes and increase endothelial cell purity 
\cite{perriere2005puromycin}. Cultures were maintained at $37\celsius$ in a humidified atmosphere of
5\% CO$_2$ / 95\% air. 
The medium was changed 24 hours after plating to remove non-adherent cells, red blood cells, and debris. At 48 hours after plating, the medium was changed again with new medium containing all the components listed above, except puromycin. The purified primary 
BMECs were used to construct in-vitro models when 80\% confluent (typically the 5th day after isolation). 

\subsection{Construction of the in-vitro blood-brain barrier model}
Monolayers of brain microvascular endothelial cells were used for all experiments. Endothelial cells were briefly 
treated with 0.25\% Trypsin-EDTA (Sigma-Aldrich) and seeded on the inside of a fibronectin-collagen type IV (0.1 mg/ml, each)
coated polyester membrane ($0.33 cm^2, 0.4 \mu m$ pore size) of a transwell-clear insert (Corning, Tewksbury MA) at a density 
of $4 \times 104$ cells per well. The medium used to plate the cells each of  the transwells fitted to a 24-well plate contained
all the components of BMEC medium, listed above, with the addition of hydrocortisone ($500 nM$; Sigma-Aldrich). The medium in 
the luminal chamber was changed 24 hours after seeding. BMEC monolayers were cultured for 3 days before use in blast experiments. 
Transendothelial electrical resistance (TEER, in $\Omega \times cm^2$) was measured using an ohmmeter equipped with an STX-2 
electrode (World Precision Instruments; Sarasota, FL). The TEER of cell-free transwell-clear inserts was subtracted from 
obtained values. TEER was measured immediately prior to blast exposure and 24 hours post-exposure.

\subsection{Exposure of BMEC to Blast} 
Transwells were placed into the blasting apparatus, consisting of a modified 24-well plate configuration containing only 
4 wells of the 24-well plates with a rubber gasket fitted to the modified plate. The medium was discarded from the luminal 
side of the transwell inserts and the inserts were placed in the middle two chambers of the blasting apparatus. The wells 
were filled completely with serum-free DMEM/F12 medium containing bFGF (1 ng/ml) and hydrocortisone (500 nM). A rubber 
gasket was placed between the filled wells and the lid of the apparatus to completely seal the chambers without air bubbles. 
The treatment apparatus (a single row of 4 transwell chambers with the two chambers in the middle containing the membrane 
inserts with BMECs) was then taped firmly to a rigid steel frame with 1/4 inch wire mesh, mounted in the blast tube, and 
exposed to a single mild blast (range: $11.0$ to $13.9$ peak $\mathrm{psi}$). Non-blasted sham controls were prepared and 
processed as above but were not exposed to a blast. Following treatment (blast or sham), the medium was aspirated from the chambers.
The inserts were placed in a 24-well plate with fresh serum-free medium and returned to $37\celsius$ in a humidified atmosphere 
of $5\%$ CO$_2/95\%$ air.

\subsection{Transendothelial permeability}
\label{sec:permeability}
Permeability to [14C]-sucrose was measured 24 hours after exposure to blast. Transwell inserts were first washed with 
physiological buffer containing $1\%$ bovine serum albumin ($141 mM$ NaCl, $4.0 mM$ KCl, $2.8 mM$
CaCl$_2$, $1.0 mM$ MgSO$_4$, 
$1.0 mM$ NaH$_2$PO$_4$, $10 mM$ HEPES, $10 mM$ D-glucose and $1\%$ BSA, pH $7.4$). The inserts were placed in a new 24-well 
plate containing $600 \mu l$ physiological buffer with $1\%$ BSA in the abluminal chamber. To initiate permeability 
experiments, [14C]-sucrose ($150,000 cpm/well$) in physiological buffer with $1\%$ BSA was added to the luminal 
chamber and $500 \mu l$ samples were collected from the abluminal chamber at 10, 20, 30, and 45 $\mathrm{min}$. When 
samples were removed from the abluminal chamber, an equal volume of fresh $1\%$ BSA/physiological buffer was immediately
added to the abluminal chamber to replace the sample volume. Liquid scintillation fluid was added to each sample and the 
radioactivity was measured using a liquid scintillation counter. The permeability coefficient and clearance of [$^{14}$C]-sucrose 
was calculated according to previously published methods \cite{dehouck1992drug}. Clearance was expressed as microliters of 
radioactive tracer diffusing from the luminal to the abluminal chamber, and it was calculated using the initial amount of radioactivity
in the loading chamber and the measured amount of radioactivity in the collected samples. Clearance ($\mu$L) = [C]C $\times$ 
VC / [C]L, where [C]L was the initial amount of radioactivity per microliter of the solution loaded into the insert 
(in $cpm/\mu L$), [C]C  was the radioactivity per microliter in the collected sample (in $cpm/ \mu l$), and VC is the volume 
of collecting chamber (in $\mu l$). The clearance volume increased linearly with time. The volume cleared was plotted versus 
time, and the slope was estimated by linear regression analysis. The slope of clearance curves for the BMEC monolayer plus 
transwell membrane was denoted by $PS_{app}$, where $PS$ is the permeability $\times$ surface area product (in $\mu L/\mathrm{min}$).
The slope of the clearance curve with a transwell membrane without BMECs was denoted by $PS_{membrane}$. The real $PS$ value for the 
BMEC monolayer ($PS_e$) was calculated from $1 / PS_{app} = 1 /PS_{membrane} + 1/PS_e$. The $PS_e$ values were divided by the 
surface area of the transwell inserts to generate the endothelial permeability coefficient ($P_e$, in $\mu l/(\mathrm{min}/ \mathrm{cm}^2)$).
Statistical analysis of TEER and sucrose permeability data was carried out using standard one-way analysis of variance (ANOVA) and were performed using SPSS software (IBM, Armonk NY). p values for correlations between blast intensity and TEER or sucrose permeability denote two-tailed statistical significance outcomes of a Pearson correlation. 

\subsection{Confocal Microscopy}
BMECs were washed in PBS and fixed with $4\%$ paraformaldehyde for 10 minutes at 4ºC. Cells were permeabilized with 
$0.1\%$ TRITON-X100 for 10 $\mathrm{min}$ at RT and blocked with $5\%$ BSA for 30 $\mathrm{min}$ at RT. They were then 
incubated for 1 hour at RT with primary antibody, ZO-1 (AbCam, Cambridge, UK) or claudin-5 (AbCam, Cambridge, UK), 
followed by incubation with Alexa Fluor 488 conjugated secondary antibody (Life Technologies, Carlsbad, CA). The monolayer-net
was then mounted on slides using Prolong Gold anti-fade with DAPI (Life Technologies, Grand Isle, NY) to stain cell nuclei. 
The monolayers were imaged using a TCS SP5 confocal microscope (Leica, Buffalo Grove, IL) with a $20\times 0.7$ numerical 
aperture objective. Only representative monolayer fields of cellular interfaces expressing claudin-5 and ZO-1 were imaged 
from 6 blast-exposed and 6-sham endothelial cultures. The monolayer-nets were imaged using a $0.2 \mu m$ z-plane step size 
for 27 slices representing a total depth of $5.4 \mu m$. Primary antibodies for claudin-5 and ZO-1 were purchased from 
Zymed (San Francisco, CA). Serial three-dimensional reconstructions of confocal images were carried out using Imaris 
software (Bitplane, South Windsor, CT). Figures were prepared using Photoshop and Imaris software using only linear 
brightness and contrast adjustments that were applied identically among control and blast-exposed specimens for each 
figure all image acquisition parameters were held constant in acquiring data for both identical control and blast-exposed 
specimens for each experiment.

\bibliographystyle{siam}
\bibliography{extrarefs.bib,library.bib,biorefs.bib}

\end{document}